\documentclass[11pt,twoside]{article}

\usepackage{amsfonts,epsfig,graphicx}
\usepackage{amsmath,amssymb,amsthm}
\usepackage{fullpage}

\usepackage{epsf}
\usepackage{epstopdf}
\usepackage{fancyheadings}
\usepackage{graphics}
\usepackage{amsfonts}
\usepackage{amsmath}
\usepackage{psfrag}
\usepackage{color}
\usepackage{algorithm,algorithmic}
\usepackage{mathtools}

\theoremstyle{plain}

\newtheorem{theo}{Theorem}[section]

\newtheorem{lem}{Lemma}[section]
\newtheorem{prop}{Proposition}[section]
\newtheorem{cor}{Corollary}[section]

\theoremstyle{definition} 

\newtheorem{nota}{Notation}[section]
\newtheorem{de}{Definition}[section]
\newtheorem{exa}{Example}[section]
\newtheorem{as}{Assumption}[section]
\newtheorem{alg}{Algorithm}[section]

\newcommand{\btheo}{\begin{theo}}
\newcommand{\bde}{\begin{de}}
\newcommand{\ble}{\begin{lem}}
\newcommand{\bpr}{\begin{prop}}
\newcommand{\bno}{\begin{nota}}
\newcommand{\bex}{\begin{exa}}
\newcommand{\bcor}{\begin{cor}}
\newcommand{\spro}{\begin{proof}}
\newcommand{\bas}{\begin{as}}
\newcommand{\balg}{\begin{alg}}

\newcommand{\etheo}{\end{theo}}
\newcommand{\ede}{\end{de}}
\newcommand{\ele}{\end{lem}}
\newcommand{\epr}{\end{prop}}
\newcommand{\eno}{\end{nota}}
\newcommand{\eex}{\end{exa}}
\newcommand{\ecor}{\end{cor}}
\newcommand{\fpro}{\end{proof}}
\newcommand{\eas}{\end{as}}
\newcommand{\ealg}{\end{alg}}

\theoremstyle{plain}

\newtheorem{theos}{Theorem}
\newtheorem{props}{Proposition}
\newtheorem{lems}{Lemma}
\newtheorem{cors}{Corollary}

\theoremstyle{definition}
\newtheorem{exas}{Example}
\newtheorem{algs}{Algorithm}
\newtheorem{asss}{Assumption}
\newtheorem{defns}{Definition}

\newcommand{\btheos}{\begin{theos}}
\newcommand{\etheos}{\end{theos}}
\newcommand{\bprops}{\begin{props}}
\newcommand{\eprops}{\end{props}}
\newcommand{\bdes}{\begin{defns}}
\newcommand{\edes}{\end{defns}}
\newcommand{\blems}{\begin{lems}}
\newcommand{\elems}{\end{lems}}
\newcommand{\bcors}{\begin{cors}}
\newcommand{\ecors}{\end{cors}}
\newcommand{\bexs}{\begin{exas}}
\newcommand{\eexs}{\end{exas}}
\newcommand{\balgs}{\begin{algs}}
\newcommand{\ealgs}{\end{algs}}
\newcommand{\bass}{\begin{asss}}
\newcommand{\eass}{\end{asss}}

\DeclarePairedDelimiter{\ceil}{\lceil}{\rceil}

\newcommand{\numobs}{\ensuremath{n}}
\newcommand{\usedim}{\ensuremath{d}}
\newcommand{\kdim}{\ensuremath{s}}

\newcommand{\xstar}{\ensuremath{x^*}}

\newcommand{\xhat}{\ensuremath{\widehat{x}}}

\newcommand{\mprob}{\ensuremath{\mathbb{P}}}

\newcommand{\numproj}{\ensuremath{m}}
\newcommand{\ConeSet}{\ensuremath{\mathcal{K}}}
\newcommand{\Sketch}{\ensuremath{S}}
\newcommand{\sketch}{\ensuremath{s}}
\newcommand{\Width}{\ensuremath{\mathcal{W}}}
\newcommand{\Constraint}{\ensuremath{\mathcal{C}}}
\newcommand{\Amat}{\ensuremath{A}}

\newcommand{\real}{\ensuremath{\mathbb{R}}}
\newcommand{\defn}{\ensuremath{: \, =}}
\newcommand{\vnew}{\ensuremath{v}}

\newcommand{\dprime}{\ensuremath{{\prime \prime}}}
\newcommand{\newdec}{\ensuremath{{\lambda}}}

\newcommand{\inprod}[2]{\ensuremath{\langle #1 , \, #2 \rangle}}
\newcommand{\sigkminsq}{\ensuremath{\gamma^-_{\kdim}}}

\newcommand{\AMIN}{\ensuremath{A_{\min}}}
\newcommand{\AMAX}{\ensuremath{A_{\max}}}
\newcommand{\cov}{\ensuremath{\operatorname{cov}}}

\newcommand{\Sphere}[1]{\ensuremath{\mathcal{S}^{#1-1}}}
\newcommand{\SPHERE}[1]{\ensuremath{\Sphere{#1}}}

\newcommand{\Term}{\ensuremath{T}}

\newcommand{\CEXP}[1]{\ensuremath{e^{#1}}}

\newcommand{\Exs}{\ensuremath{\mathbb{E}}}

\newcommand{\diag}{\ensuremath{\mbox{diag}}}

\newcommand{\fixvec}{\ensuremath{u}}

\newcommand{\ENCMIN}[1]{\ensuremath{\Big \{#1 \Big \}}}

\newcommand{\xit}[1]{\ensuremath{x^{#1}}}
\newcommand{\SketchIt}[1]{\ensuremath{\Sketch^{#1}}}
\newcommand{\ZSUP}{\ensuremath{Z_2}}
\newcommand{\ZINF}{\ensuremath{Z_1}}

\newlength{\widebarargwidth}
\newlength{\widebarargheight}
\newlength{\widebarargdepth}

\newcommand{\KCONE}{\ensuremath{\mathcal{K}}}

\newcommand{\Event}{\ensuremath{\mathcal{E}}}

\newcommand{\widgraph}[2]{\includegraphics[keepaspectratio,width=#1]{#2}}

\newcommand{\hessh}{\ensuremath{\nabla^2 f(x^t)^{1/2}}}

\newcommand{\dom}{\ensuremath{{\rm dom\,}}}

\newcommand{\fp}{\ensuremath{\nabla f}}
\newcommand{\fpp}{\ensuremath{\nabla^2 f}}

\newcommand{\order}{\ensuremath{\mathcal{O}}}

\DeclareMathOperator{\trace}{trace}
\DeclareMathOperator{\rank}{rank}

\theoremstyle{definition}
\newtheorem{example}{Example}

\newcommand{\HessSqrt}[1]{\ensuremath{\nabla^2 f(#1)^{1/2}}}
\newcommand{\HessSqrtBase}[1]{\ensuremath{\nabla^2 \fbase(#1)^{1/2}}}

\newcommand{\HessSqrtT}[1]{\ensuremath{(\nabla^2 f(#1)^{1/2})^T}}
\newcommand{\Hess}[1]{\ensuremath{\nabla^2 f(#1)}}
\newcommand{\xittil}[1]{\ensuremath{\tilde{x}^{#1}}}

\newcommand{\myvec}{\ensuremath{\mbox{vec}}}

\newcommand{\SPECFUN}[1]{\ensuremath{\tilde \Phi(#1)}}
\newcommand{\SPECFUNTWO}[2]{\ensuremath{\Phi(#1; #2)}}

\newcommand{\relow}{\ensuremath{\gamma}}
\newcommand{\reup}{\ensuremath{\beta}}

\newcommand{\LipCon}{\ensuremath{L}}

\newcommand{\TotalIt}{\ensuremath{N}}

\newcommand{\linea}{\ensuremath{a}}
\newcommand{\lineb}{\ensuremath{b}}

\newcommand{\HackCone}[1]{\HessSqrt{#1} \ConeSet}

\newcommand{\DelIt}[1]{\ensuremath{\Delta^{#1}}}

\newcommand{\CONNUM}{\ensuremath{r}}
\newcommand{\NEWGAM}{\ensuremath{\nu}}
\newcommand{\NEWETA}{\ensuremath{\eta}}

\newcommand{\ConeSetIt}[1]{\ensuremath{\ConeSet^{#1}}}
\newcommand{\numprojit}[1]{\ensuremath{\numproj^{#1}}}

\newcommand{\WidthSq}{\ensuremath{\Width^2}}

\newcommand{\matsnorm}[2]{|\!|\!| #1 | \! | \!|_{{#2}}}

\newcommand{\opnorm}[1]{\ensuremath{\matsnorm{#1}{\tiny{\mbox{op}}}}}

\newcommand{\Zone}{\ensuremath{Z_1}}
\newcommand{\Ztwo}{\ensuremath{Z_2}}

\newcommand{\ZONEHACK}[2]{\ensuremath{Z_1(#1 ; \, #2)}}
\newcommand{\ZTWOHACK}[2]{\ensuremath{Z_2(#1 ; \, #2)}}
\newcommand{\ZONET}{\ensuremath{Z_1^t}}
\newcommand{\ZTWOT}{\ensuremath{Z_2^t}}

\newcommand{\RHS}{\ensuremath{\mbox{RHS}}}
\newcommand{\vnewt}{\ensuremath{v_{\tiny{\mbox{NE}}}}}
\newcommand{\vsketch}{\ensuremath{v_{\tiny{\mbox{NSK}}}}}

\newcommand{\xnewt}{\ensuremath{x_{\tiny{\mbox{NE}}}}}
\newcommand{\xsketch}{\ensuremath{x_{\tiny{\mbox{NSK}}}}}

\newcommand{\newdecsketch}{\ensuremath{\widetilde{\lambda}}}

\newcommand{\SHORTA}{\newdecsketch}

\newcommand{\stepsize}{\ensuremath{s}}
\newcommand{\Devent}{\ensuremath{\mathcal{D}}}

\newcommand{\gfun}{\ensuremath{g}}
\newcommand{\shat}{\ensuremath{\widehat{u}}}
\newcommand{\ustep}{\ensuremath{u}} 

\newcommand{\errvec}{\ensuremath{\widehat{e}}}

\newcommand{\rdim}{\ensuremath{r}}

\newcommand{\fbase}{\ensuremath{f_0}}

\newcommand{\Lset}{\ensuremath{\mathcal{L}}}

\newcommand{\Yspace}{\ensuremath{\mathcal{Y}}}
\newcommand{\LinkFun}{\ensuremath{\psi}}
\newcommand{\LinkFunDouble}{\ensuremath{\LinkFun^{\prime \prime}}}
\newcommand{\LinkFunSingle}{\ensuremath{\LinkFun^{\prime}}}

\newcommand{\fzero}{\ensuremath{f_0}}
\newcommand{\gcon}{\ensuremath{g}}

 \long\def\comment#1{}

\makeatletter
\long\def\@makecaption#1#2{
        \vskip 0.8ex
        \setbox\@tempboxa\hbox{\small {\bf #1:} #2}
        \parindent 1.5em  
        \dimen0=\hsize
        \advance\dimen0 by -3em
        \ifdim \wd\@tempboxa >\dimen0
                \hbox to \hsize{
                        \parindent 0em
                        \hfil 
                        \parbox{\dimen0}{\def\baselinestretch{0.96}\small
                                {\bf #1.} #2
                                } 
                        \hfil}
        \else \hbox to \hsize{\hfil \box\@tempboxa \hfil}
        \fi
        }
\makeatother

\begin{document}

\begin{center} 

{\LARGE{\bf{ Newton Sketch: A Linear-time Optimization Algorithm with
      Linear-Quadratic Convergence}}} \\

  \vspace{1cm}
\begin{tabular}{ccc}
  {\large Mert Pilanci$^{1}$} & & {\large{Martin J.\ Wainwright$^{1,2}$}}
\end{tabular}

 \vspace{.2cm}
  \texttt{\{mert, wainwrig\}@berkeley.edu} \\
  \vspace{1cm}
  {\large University of California, Berkeley} \\
  \vspace{.15cm} $^1$Department of Electrical Engineering and Computer
  Science ~~~~ $^2$Department of Statistics

\vspace*{.2in}

\today

\end{center}

\vspace*{.5in}


\begin{abstract}
We propose a randomized second-order method for optimization known as
the Newton Sketch: it is based on performing an approximate Newton
step using a randomly projected or sub-sampled Hessian. For
self-concordant functions, we prove that the algorithm has
super-linear convergence with exponentially high probability, with
convergence and complexity guarantees that are independent of
condition numbers and related problem-dependent quantities.  Given a
suitable initialization, similar guarantees also hold for strongly
convex and smooth objectives without self-concordance.  When
implemented using randomized projections based on a sub-sampled
Hadamard basis, the algorithm typically has substantially lower complexity than
Newton's method. We also describe extensions of our methods to
programs involving convex constraints that are equipped with
self-concordant barriers.  We discuss and illustrate applications to
linear programs, quadratic programs with convex constraints, logistic
regression and other generalized linear models, as well as
semidefinite programs.
\end{abstract}


\section{Introduction}

Relative to first-order methods, second-order methods for convex
optimization enjoy superior convergence in both theory and practice.
For instance, Newton's method converges at a quadratic rate for
strongly convex and smooth problems, and moreover, even for weakly
convex functions (i.e. not strongly convex), modifications of Newton's method has super-linear convergence compared to
the much slower $1/T^2$ convergence rate that can be achieved by a
first-order method like accelerated gradient descent (see e.g. \cite{yamashita2001rate}).  More
importantly, at least in a uniform sense, the $1/T^2$-rate is known to
be unimprovable for first-order methods~\cite{Nesterov04}. Yet another
issue in first-order methods is the tuning of step size, whose optimal
choice depends on the strong convexity parameter and/or smoothness of
the underlying problem.  For example, consider the problem of
optimizing a function of the form $x \mapsto g(\Amat x)$, where $\Amat
\in \real^{\numobs \times \usedim}$ is a ``data matrix'', and $g:
\real^\numobs \rightarrow \real$ is a twice-differentiable function.
Here the performance of first-order methods will depend on both the
convexity/smoothness of $g$, as well as the conditioning of the data
matrix.  In contrast, whenever the function $g$ is self-concordant,
then Newton's method with suitably damped steps has a global complexity guarantee that is
provably independent of such problem-dependent parameters.

On the other hand, each step of Newton's method requires solving a
linear system defined by the Hessian matrix. For instance, in
application to the problem family just described involving an $\numobs
\times \usedim$ data matrix, each of these steps has complexity
scaling as $\order( \numobs \usedim^2)$.  For this reason, both
forming the Hessian and solving the corresponding linear system pose a
tremendous numerical challenge for large values of $(\numobs,
\usedim)$--- for instance, values of thousands to millions, as is
common in big data applications, In order to address this issue, a
multitude of different approximations to Newton's method have been
proposed and studied in the literature.  Quasi-Newton methods form
estimates of the Hessian by successive evaluations of the gradient
vectors and are computationally cheaper. Examples of such methods
include DFP and BFGS schemes and also their limited memory versions
(see the book~\cite{wright1999numerical} for further details). A
disadvantage of such approximations based on first-order information
is that the associated convergence guarantees are typically much
weaker than those of Newton's method and require stronger assumptions.
Under restrictions on the eigenvalues of the Hessian (strong convexity
and smoothness), Quasi-Newton methods typically exhibit local
super-linear convergence.

In this paper, we propose and analyze a randomized approximation of
Newton's method, known as the \emph{Newton Sketch}.  Instead of
explicitly computing the Hessian, the Newton Sketch method
approximates it via a random projection of dimension $\numproj$.  When
these projections are carried out using the randomized Hadamard
transform, each iteration has complexity $\order(\numobs \usedim
\log(\numproj) + \usedim \numproj^2)$. Our results show that it is
always sufficient to choose $\numproj$ proportional to $\min \{
\usedim, \numobs \}$, and moreover, that the sketch dimension $m$ can be
much smaller for certain types of constrained problems.  Thus, in the
regime $\numobs > \usedim$ and with $\numproj \asymp \usedim$, the
complexity per iteration can be substantially lower than the
$\order(\numobs \usedim^2)$ complexity of each Newton step. Specifically for $n\ge d^2$, the complexity of Newton Sketch per iteration is $\order(nd\log d)$, which is linear in the input size ($nd$) and comparable to first order methods which only access the derivative $g^\prime(Ax)$. Moreover, we show that for self-concordant functions, the total complexity of obtaining a $\delta$-optimal solution is $\order(nd\log d \log(1/\delta))$, and does not depend on constants such as strong convexity or smoothness parameters unlike first order methods. On the
other hand, for problems with $\usedim > \numobs$, we also provide a
dual strategy which effectively has the same guarantees with roles of
$\usedim$ and $\numobs$ exchanged.

We also consider other random projection matrices and sub-sampling
strategies, including partial forms of random projection that exploit
known structure in the Hessian. For self-concordant functions, we
provide an affine invariant analysis proving that the convergence is
linear-quadratic and the guarantees are independent of the function
and data, such as condition numbers of matrices involved in the
objective function.  Finally, we describe an interior point method to
deal with arbitrary convex constraints which combines the Newton
sketch with the barrier method. We provide an upper bound on the total
number of iterations required to obtain a solution with a
pre-specified target accuracy.

The remainder of this paper is organized as follows.  We begin in
Section~\ref{SecBackground} with some background on the classical form
of Newton's method, random matrices for sketching, and Gaussian widths
as a measure of the size of a set.  In Section~\ref{SecNewtonSketch},
we formally introduce the Newton Sketch, including both fully and
partially sketched versions for unconstrained and constrained
problems.  We provide some illustrative examples in
Section~\ref{SecExamples} before turning to local convergence theory
in Section~\ref{SecLocal}.  Section~\ref{SecGlobal} is devoted to
global convergence results for self-concordant functions, in both the
constrained and unconstrained settings.  In
Section~\ref{SecApplications}, we consider a number of applications
and provide additional numerical results.  The bulk of our proofs are
in given in Section~\ref{SecProofs}, with some more technical aspects
deferred to the appendices.


\section{Background}
\label{SecBackground}

We begin with some background material on the standard form of
Newton's method, various types of random sketches, and the notion of
Gaussian width as a complexity measure.

\subsection{Classical version of Newton's method}
\label{SecBackNewton}

In this section, we briefly review the convergence properties and
complexity of the classical form of Newton's method; see the
sources~\cite{wright1999numerical,Boyd02,Nesterov04} for further background.

Let $f: \real^\usedim \rightarrow \real$ be a closed,
convex and twice-differentiable function that is bounded below.  Given
a convex set $\Constraint$, we assume that the constrained minimizer
\begin{align}
\xstar & \defn \arg \min_{x \in \Constraint} f(x)
\end{align}
is uniquely defined, and we define the minimum and maximum eigenvalues
$\relow = \lambda_{min}(\nabla^2 f(\xstar))$ and $\reup =
\lambda_{max}(\nabla^2 f(\xstar))$ of the Hessian evaluated at the
minimum.

We assume moreover that the Hessian map $x \mapsto \nabla^2 f(x)$ is
Lipschitz continuous with modulus $\LipCon$, meaning that
\begin{align}
\label{EqnDefnHessianLipCon}
\opnorm{\Hess{x + \Delta} - \Hess{x}} & \leq \LipCon \, \|\Delta\|_2.
\end{align}
Under these conditions and given an initial point $\xittil{0} \in
\Constraint$ such that $\|\xittil{0} - \xstar\|_2 \leq \frac{\relow}{2
  \LipCon}$, the Newton updates are guaranteed to converge
quadratically---viz.
\begin{align*}
\|\xittil{t+1} - \xstar\|_2 & \leq \frac{2 \LipCon}{\relow}
  \|\xittil{t} - \xstar\|_2^2,
\end{align*}
This result is classical: for instance, see Boyd and
Vandenberghe~\cite{Boyd02} for a proof. Newton's method can be
slightly modified to be globally convergent by choosing the step sizes
via a simple backtracking line-search procedure.

The following result characterizes the complexity of Newton's method
when applied to self-concordant functions and is central in the
development of interior point methods (for instance, see the
books~\cite{NesNem94,Boyd02}). We defer the definitions of
self-concordance and the line-search procedure in the following
sections. The number of iterations needed to obtain a $\delta$
approximate minimizer of a strictly convex self-concordant function
$f$ is bounded by
\begin{align*}
\frac{20 - 8 \linea}{\linea \lineb (1 - 2 \linea)} \left( f(\xit{0}) -
f(\xstar) \right) + \log_2 \log_2(1/\delta)\,,
\end{align*}
where $\linea, \lineb$ are constants in the line-search
procedure.\footnote {Typical values of these constants are $\linea =
  0.1$ and $\lineb = 0.5$.}


\subsection{Different types of randomized sketches}
\label{SecBackSketches}

Various types of randomized sketches are possible, and we describe a
few of them here.  Given a sketching matrix $\Sketch \in
\real^{\numproj \times \numobs}$, we use
$\{\sketch_i\}_{i=1}^\numproj$ to denote the collection of its
$\numobs$-dimensional rows.  We restrict our attention to sketch
matrices that are zero-mean, and that are normalized so that
$\Exs[\Sketch^T \Sketch/\numproj] = I_\numobs$.

\paragraph{Sub-Gaussian sketches:}  The most classical sketch
is based on a random matrix \mbox{$\Sketch \in \real^{\numproj \times
    \numobs}$} with i.i.d. standard Gaussian entries, or somewhat more
generally, sketch matrices based on i.i.d. sub-Gaussian rows.  In
particular, a zero-mean random vector $\sketch \in \real^{\numobs}$ is
$1$-sub-Gaussian if for any $u \in \real^\numobs$, we have
\begin{align}
\mprob[ \inprod{\sketch}{u} \geq \epsilon \|u\|_2 \big] & \leq
\CEXP{-\epsilon^2/2} \qquad \mbox{for all $\epsilon \geq 0$.}
\end{align}
For instance, a vector with i.i.d. $N(0,1)$ entries is
$1$-sub-Gaussian, as is a vector with i.i.d. Rademacher entries
(uniformly distributed over $\{-1, +1\}$).  We use the terminology
\emph{sub-Gaussian sketch} to mean a random matrix $\Sketch \in
\real^{\numproj \times \numobs}$ with i.i.d.  rows that are zero-mean,
$1$-sub-Gaussian, and with $\cov(\sketch) = I_\numobs$.

From a theoretical perspective, sub-Gaussian sketches are attractive
because of the well-known concentration properties of sub-Gaussian
random matrices (e.g.,~\cite{DavSza01,Ver11}).  On the other hand,
from a computational perspective, a disadvantage of sub-Gaussian
sketches is that they require matrix-vector multiplications with
unstructured random matrices.  In particular, given a data matrix
$\Amat \in \real^{\numobs \times \usedim}$, computing its sketched
version $\Sketch \Amat$ requires $\order(\numproj \numobs \usedim)$
basic operations in general (using classical matrix multiplication).


\paragraph{Sketches based on randomized orthonormal systems (ROS):}

The second type of randomized sketch we consider is \emph{randomized
  orthonormal system} (ROS), for which matrix multiplication can be
performed much more efficiently.  In order to define a ROS sketch, we
first let $H \in \real^{\numobs \times \numobs}$ be an orthonormal
matrix with entries $H_{ij} \in [ -\frac{1}{\sqrt{\numobs}},
  \frac{1}{\sqrt{\numobs}} ]$.  Standard classes of such matrices are
the Hadamard or Fourier bases, for which matrix-vector multiplication
can be performed in $\order(\numobs \log \numobs)$ time via the fast
Hadamard or Fourier transforms, respectively.  Based on any such
matrix, a sketching matrix $\Sketch \in \real^{\numproj \times
  \numobs}$ from a ROS ensemble is obtained by sampling i.i.d. rows of
the form
\begin{align*}
\sketch^T & = \sqrt{\numobs} e_j^T H D \qquad \mbox{with probability
  $1/\numobs$ for $j = 1, \ldots, \numobs$},
\end{align*}
where the random vector $e_j \in \real^\numobs$ is chosen uniformly at
random from the set of all $\numobs$ canonical basis vectors, and $D =
\diag(\nu)$ is a diagonal matrix of i.i.d. Rademacher variables
\mbox{$\nu \in \{-1, +1\}^\numobs$.}  Given a fast routine for
matrix-vector multiplication, the sketch $\Sketch M$ for a data matrix
$M\in \real^{\numobs \times \usedim}$ can be formed in $\order(\numobs
\, \usedim \log \numproj)$ time (for instance, see the
papers~\cite{Ailon08,ailon2006approximate}).


\paragraph{Sketches based on random row sampling:}
Given a probability distribution $\{p_j\}_{j=1}^\numobs$ over
$[\numobs] = \{1, \ldots, \numobs \}$, another choice of sketch is to
randomly sample the rows of a data matrix $M$ a total of $\numproj$
times with replacement from the given probability distribution.  Thus,
the rows of $S$ are independent and take on the values
\begin{align*}
s^T & = \frac{e_j}{\sqrt{p_j}} \qquad \mbox{with probability $p_j$ for
  $j = 1, \ldots, \numobs$}
\end{align*}
where $e_j \in \real^\numobs$ is the $j^{th}$ canonical basis vector.
Different choices of the weights $\{p_j\}_{j=1}^\numobs$ are possible,
including those based on the row $\ell_2$ norms $p_j \propto
\|Me_j\|_2^2$ and leverage values of $M$---i.e., $p_j \propto
\|Ue_j\|_2$ for $j= 1, \ldots, \numobs$, where $U \in \real^{\numobs
  \times \usedim}$ is the matrix of left singular vectors of $M$
\cite{drineas2012fast}. When $M\in \real^{n\times d}$ is the adjacency
matrix of a graph with $d$ vertices and $n$ edges, the leverage scores
of $M$ are also known as effective resistances which can be used to
sub-sample edges of a given graph by preserving its spectral
properties~\cite{spielman2011graph}.


\subsection{Gaussian widths}
\label{SecBackWidths}

In this section, we introduce some background on the notion of
Gaussian width, a way of measuring the size of a compact set in
$\real^\usedim$.  These width measures play a key role in the analysis
of randomized sketches.  Given a compact subset $\Lset \subseteq
\real^d$, its Gaussian width is given by
\begin{align}
\label{EqnDefnGaussWidth}
\Width(\Lset) & \defn \Exs_g \big[ \max_{z \in \Lset} |
  \inprod{g}{z} |\big]
\end{align}
where $g \in \real^\numobs$ is an i.i.d. sequence of $N(0,1)$
variables.  This complexity measure plays an important role in Banach
space theory, learning theory and statistics
(e.g.,~\cite{Pisier86,LedTal91,Bar05}).

Of particular interest in this paper are sets $\Lset$ that are
obtained by intersecting a given cone $\KCONE$ with the Euclidean
sphere $\SPHERE{\usedim} = \{ z \in \real^\numobs \, \mid \, \|z\|_2 =
1 \}$.  It is easy to show that the Gaussian width of any such set is
at most $\sqrt{\usedim}$, but the it can be substantially smaller,
depending on the nature of the underlying cone.  For instance, if
$\KCONE$ is a subspace of dimension $\rdim < \usedim$, then a simple
calculation yields that $\Width(\KCONE \cap \Sphere{\usedim}) \leq
\sqrt{\rdim}$.


\section{Newton sketch and local convergence}
\label{SecNewtonSketch}

With the basic background in place, let us now introduce the Newton
sketch algorithm, and then develop a number of convergence guarantees
associated with it.  It applies to an optimization problem of the form
$\min_{x \in \Constraint} f(x)$, where $f:\real^\usedim \rightarrow
\real$ is a twice-differentiable convex function, and $\Constraint
\subseteq \real^\usedim$ is a convex constraint set.


\subsection{Newton sketch algorithm}

In order to motivate the Newton sketch algorithm, recall the standard
form of Newton's algorithm: given a current iterate $\xittil{t} \in
\Constraint$, it generates the new iterate $\xittil{t+1}$ by
performing a constrained minimization of the second order Taylor
expansion---viz.
\begin{subequations}
\begin{align}
\label{EqnNewton}
\xittil{t+1} & = \arg \min_{x \in \Constraint} \ENCMIN{\frac{1}{2}
  \inprod{x - \xittil{t}}{ \Hess{\xittil{t}} \, (x - \xittil{t})} +
  \inprod{\nabla f(\xittil{t})}{x - \xittil{t}}}.
\end{align}
In the unconstrained case---that is, when $\Constraint =
\real^\usedim$---it takes the simpler form
\begin{align}
\xittil{t+1} & = \xittil{t} - \big[ \Hess{\xittil{t}} \big]^{-1}
\nabla f(\xittil{t})\,.
\end{align}
\end{subequations}

Now suppose that we have available a Hessian matrix square root
$\HessSqrt{x}$---that is, a matrix $\HessSqrt{x}$ of dimensions
$\numobs \times \usedim$ such that
\begin{align*}
\HessSqrtT{x} \HessSqrt{x} = \Hess{x} \qquad \mbox{for some integer
  $\numobs \geq \rank(\Hess{x})$.}
\end{align*}
In many cases, such a matrix square root can be computed efficiently.
For instance, consider a function of the form $f(x) = g(\Amat x)$
where $\Amat \in \real^{\numobs \times \usedim}$, and the function $g:
\real^\numobs \rightarrow \real$ has the separable form $g(\Amat x) =
\sum_{i=1}^\numobs g_i(\inprod{a_i}{x})$.  In this case, a suitable
Hessian matrix square root is given by the $\numobs \times \usedim$
matrix $\HessSqrt{x} \defn \diag \big\{ g_i^{\prime
  \prime}(\inprod{a_i}{x}) \big \}_{i=1}^\numobs \Amat$.  In
Section~\ref{SecExamples}, we discuss various concrete instantiations
of such functions.

In terms of this notation, the ordinary Newton update can be
re-written as
\begin{align*}
\xittil{t+1} & = \arg \min_{x \in \Constraint} \ENCMIN{
  \underbrace{\frac{1}{2} \|\HessSqrt{\xittil{t}} (x - \xittil{t})
    \|_2^2 + \inprod{\nabla f(\xittil{t})}{x -
      \xittil{t}}}_{\SPECFUN{x}} },
\end{align*}
and the Newton Sketch algorithm is most easily understood based on
this form of the updates.  More precisely, for a sketch dimension
$\numproj$ to be chosen, let $\Sketch \in \real^{\numproj \times
  \numobs}$ be an isotropic sketch matrix, satisfying the relation
$\Exs[\Sketch^T \Sketch] = I_{\numobs}$.  The \emph{Newton Sketch
  algorithm} generates a sequence of iterates
$\{\xit{t}\}_{t=0}^\infty$ according to the recursion
\begin{align}
\label{EqnNewtonSketch}
\xit{t+1} & \defn \arg \min_{x \in \Constraint}
\ENCMIN{\underbrace{\frac{1}{2}\| \SketchIt{t} \HessSqrt{\xit{t}} (x-
    \xit{t})\|_2^2 + \inprod{\nabla f(\xit{t})}{x -
      \xit{t}}}_{\SPECFUNTWO{x}{\SketchIt{t}}}},
\end{align}
where $\SketchIt{t} \in \real^{\numproj \times \usedim}$ is an
independent realization of a sketching matrix.  When the problem is
unconstrained, i.e., $\Constraint = \real^\usedim$ and the matrix
$\hessh (\SketchIt{t})^T \SketchIt{t} \hessh$ is invertible, the
Newton sketch update takes the simpler form to
\begin{align}
\label{EqnNewtonSketchUnconstrained}
x^{t+1 } & = \xit{t} - \left ( \hessh (\SketchIt{t})^T \SketchIt{t} \hessh
\right )^{-1} \nabla f(\xit{t}).
\end{align}
The intuition underlying the Newton sketch updates is as follows: the
iterate $\xit{t+1}$ corresponds to the constrained minimizer of the
random objective function $\SPECFUNTWO{x}{\SketchIt{t}}$ whose
expectation $\Exs[\SPECFUNTWO{x}{\SketchIt{t}}]$, taking averages over
the isotropic sketch matrix $\SketchIt{t}$, is equal to the original
Newton objective $\SPECFUN{x}$.  Consequently, it can be seen as a
stochastic form of the Newton update.

In this paper, we also analyze a \emph{partially sketched Newton
  update}, which takes the following form.  Given an additive
decomposition of the form $f = \fbase + g$, we perform a sketch of of
the Hessian $\nabla^2 \fbase$ while retaining the exact form of the
Hessian $\nabla^2 g$.  This leads to the partially sketched update
\begin{align}
\label{EqnNewtonSketchPartial}
\xit{t+1} & \defn \arg \min_{x \in \Constraint} \ENCMIN{\frac{1}{2} (x
  - \xit{t})^T Q^t (x - \xit{t}) + \inprod{\nabla f(\xit{t})}{x -
    \xit{t}}}
\end{align}
where $Q^t \defn (\SketchIt{t} \HessSqrtBase{\xit{t}})^T \SketchIt{t}
\HessSqrtBase{\xit{t}} + \nabla^2 g(\xit{t})$.

For either the fully sketched~\eqref{EqnNewtonSketch} or partially
sketched updates~\eqref{EqnNewtonSketchPartial}, our analysis shows
that there are many settings in which the sketch dimension $\numproj$
can be chosen to be substantially smaller than $\numobs$, in which
cases the sketched Newton updates will be much cheaper than a standard
Newton update.  For instance, the unconstrained
update~\eqref{EqnNewtonSketchUnconstrained} can be computed in at most
$\order(\numproj \usedim^2)$ time, as opposed to the $\order(\numobs
\usedim^2)$ time of the standard Newton update.  In constrained
settings, we show that the sketch dimension $\numproj$ can often be
chosen even smaller---even $\numproj \ll \usedim$---which leads to
further savings.

\subsection{Some examples}
\label{SecExamples}

In order to provide some intuition, let us provide some simple
examples to which the sketched Newton updates can be applied.

\begin{example}[Newton sketch for LP solving]
\label{ExaLPSolve}
Consider a linear program (LP) in the standard form
\begin{align}
\label{EqnPolytopeConstrained}
\min_{\Amat x \leq b} \, \inprod{c}{x}
\end{align}
where $\Amat \in \real^{\numobs \times \usedim}$ is a given constraint
matrix.  We assume that the polytope $\{x \in \real^\usedim \, \mid A
x \leq b \}$ is bounded so that the minimum achieved.  A barrier
method approach to this LP is based on solving a sequence of problems
of the form
\begin{align*}
\min_{x \in \real^\usedim} \ENCMIN{ \underbrace{\tau \, \inprod{c}{x}
    - \sum_{i=1}^\numobs \log(b_i- \inprod{a_i}{x})}_{f(x)}},
\end{align*}
where $a_i \in \real^\usedim$ denotes the $i^{th}$ row of $\Amat$, and
$\tau > 0$ is a weight parameter that is adjusted during the
algorithm.  By inspection, the function $f: \real^\usedim \rightarrow
\real \cup \{+\infty \}$ is twice-differentiable, and its Hessian is
given by $\Hess{x} = \Amat^T \diag \big \{ \frac{1}{(b_i -
  \inprod{a_i}{x} )^2} \big \} \Amat$.  A Hessian square root is given
by $\HessSqrt{x} \defn \diag \left(\frac{1}{|b_i- \inprod{a_i}{x}|}
\right) \Amat$, which allows us to compute a sketched version of the
Hessian square root
\begin{align*}
\Sketch \HessSqrt{x} & = \Sketch \; \diag \left(\frac{1}{|b_i -
  \inprod{a_i}{x}|} \right) \Amat.
\end{align*}
With a ROS sketch matrix, computing this matrix requires $\order
(\numobs \usedim \log(\numproj))$ basic operations.  The complexity of
each Newton sketch iteration scales as $\order(\numproj \usedim^2)$,
where $\numproj$ is at most $\usedim$. In contrast, the standard
unsketched form of the Newton update has complexity $\order(\numobs
\usedim^2)$, so that the sketched method is computationally cheaper whenever there are
more constraints than dimensions ($\numobs > \usedim$).

By increasing the barrier parameter $\tau$, we obtain a sequence of
solutions that approach the optimum to the LP, which we refer to as
the central path.  As a simple illustration,
Figure~\ref{FigCentralPath} compares the central paths generated by
the ordinary and sketched Newton updates for a polytope defined by
$\numobs = 32$ constraints in dimension $\usedim = 2$.  Each row shows
three independent trials of the method for a given sketch dimension
$\numproj$; the top, middle and bottom rows correspond to sketch
dimensions $\numproj \in \{\usedim, 4 \usedim, 16 \usedim\}$
respectively.  Note that as the sketch dimension $\numproj$ is
increased, the central path taken by the sketched updates converges to
the standard central path.

\newcommand{\mysize}{.7\textwidth}
\begin{figure}[h!]
\begin{center}
\begin{tabular}{c}
\includegraphics[width=\mysize,trim=25mm 56mm 15mm
  24mm,clip]{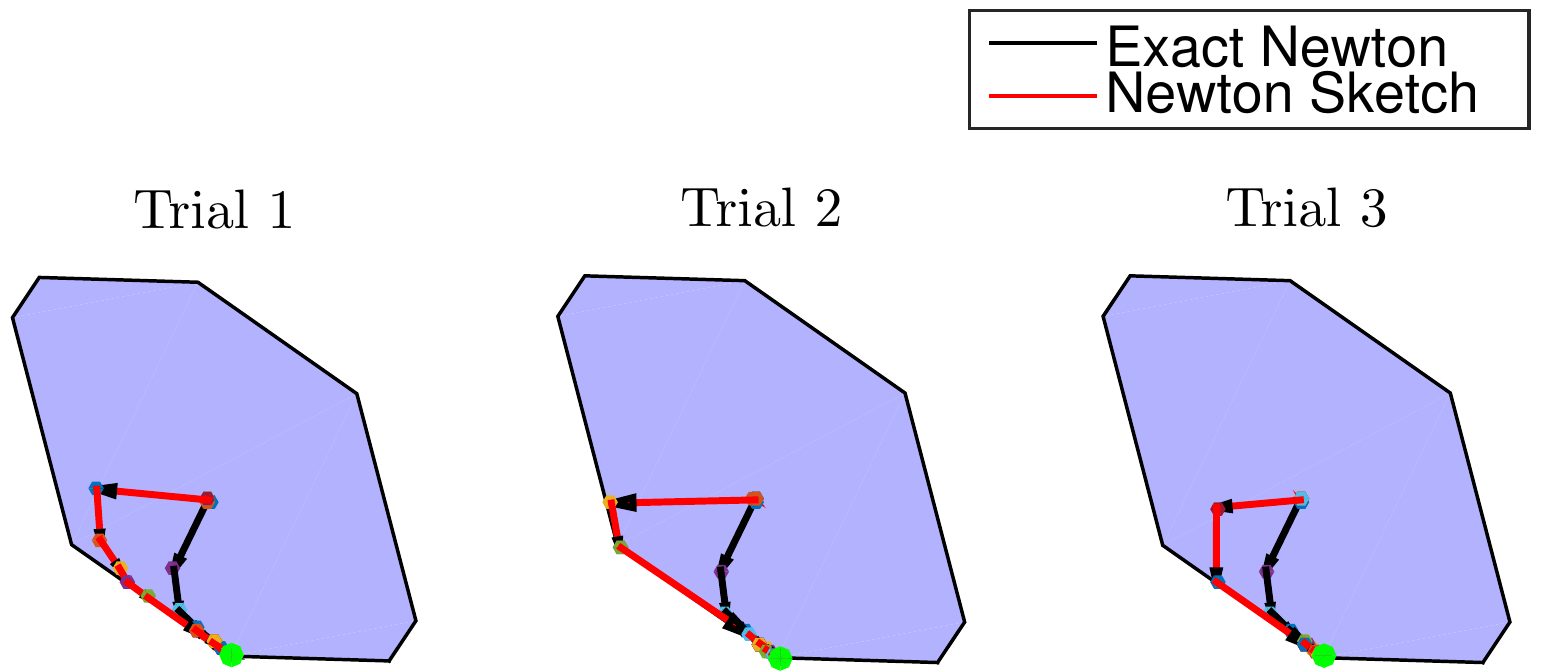} \\ (a) sketch size $\numproj =
\usedim$\\
\includegraphics[width=\mysize,trim=25mm 56mm 15mm
  24mm,clip]{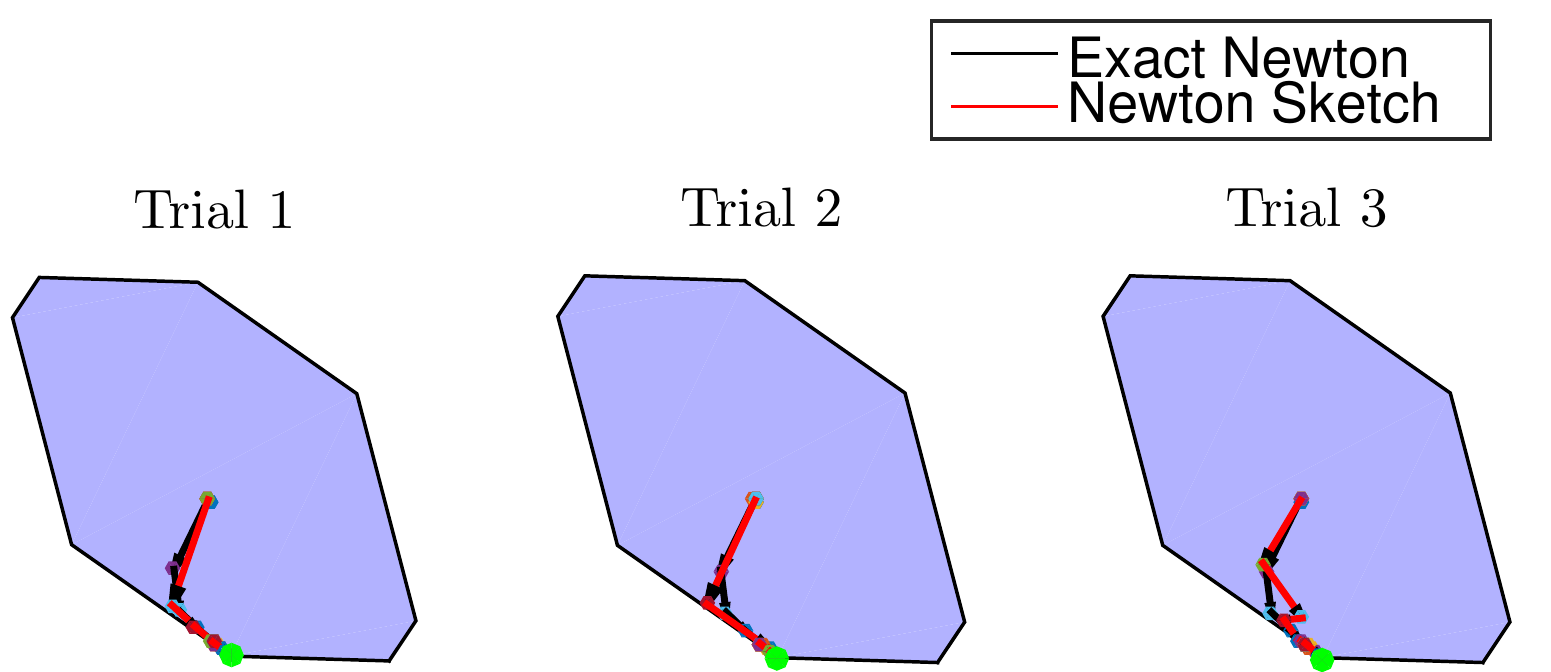} \\ (b) sketch size $\numproj =
4 \usedim$ \\
\includegraphics[width=\mysize,trim=25mm 56mm 15mm
  24mm,clip]{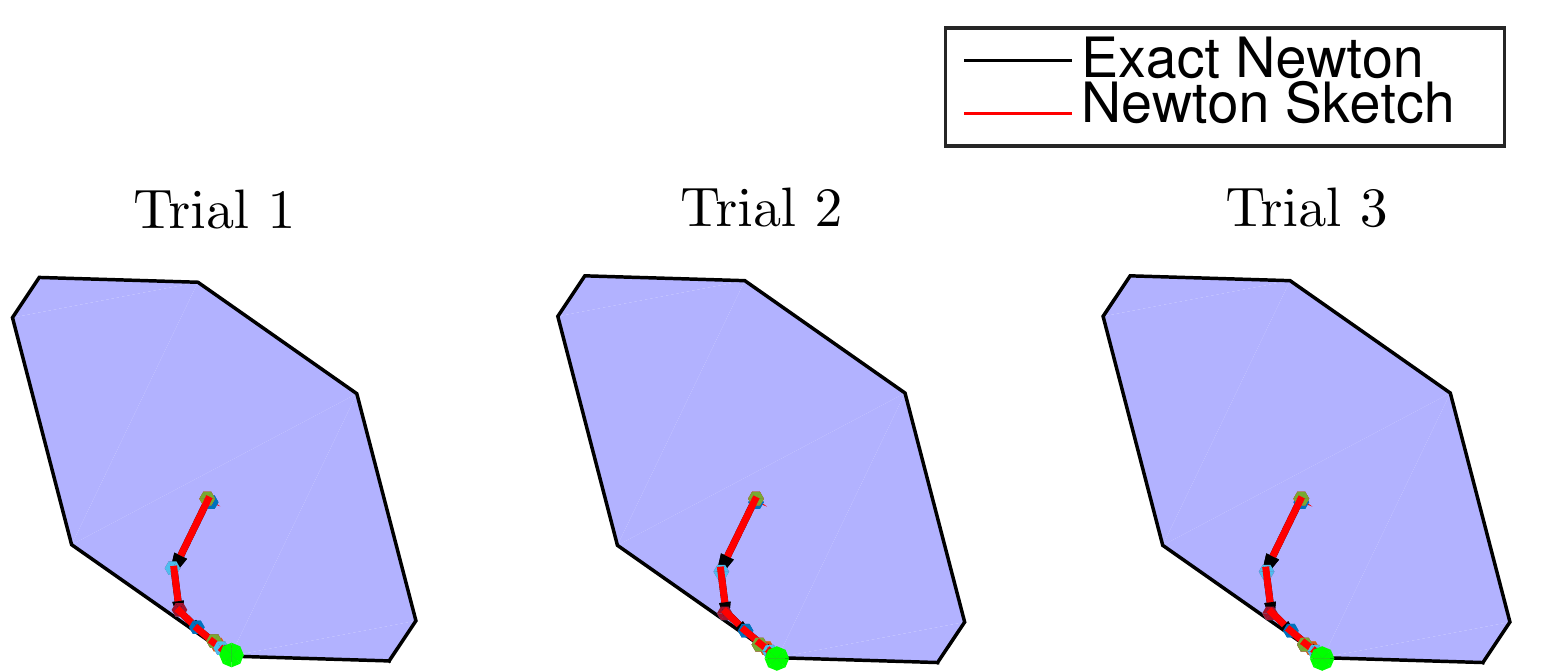} \\ (c) sketch size $\numproj
= 16 \usedim$
\end{tabular}

\caption{ Comparisons of central paths for a simple linear program in
  two dimensions.  Each row shows three independent trials for a given
  sketch dimension: across the rows, the sketch dimension ranges as
  $\numproj \in \{ \usedim, 4 \usedim, 16 \usedim\}$. The black arrows
  show Newton steps taken by the standard interior point method,
  whereas red arrows show the steps taken by the sketched version.
  The green point at the vertex represents the optimum. In all cases,
  the sketched algorithm converges to the optimum, and as the sketch
  dimension $\numproj$ increases, the sketched central path converges
  to the standard central path.}
\label{FigCentralPath}
\end{center}
\end{figure}
\end{example}

As a second example, we consider the problem of maximum likelihood
estimation for generalized linear models.

\begin{example}[Newton sketch for maximum likelihood estimation]
\label{ExaGLMs}
The class of generalized linear models (GLMs) is used to model a wide
variety of prediction and classification problems, in which the goal
is to predict some output variable $y \in \Yspace$ on the basis of a
covariate vector $a \in \real^\usedim$.  it includes as special cases
the standard linear Gaussian model (in which $\Yspace = \real$), as
well as logistic models for classification (in which $\Yspace = \{-1,
+1\}$), as well as as Poisson models for count-valued responses (in
which $\Yspace = \{0,1,2, \ldots \}$).  See the book~\cite{McCullagh}
for further details and applications.

Given a collection of $\numobs$ observations $\{(y_i,
a_i)\}_{i=1}^\numobs$ of response-covariate pairs from some GLM, the
problem of constrained maximum likelihood estimation be written in the
form
\begin{align}
\label{EqnMLEGLM}
\min_{x \in \Constraint}~ \Big \{ \underbrace{ \sum_{i=1}^\numobs
  \LinkFun(\inprod{a_i}{x}, y_i)} \Big \},
\end{align}
where $\LinkFun: \real \times \Yspace \rightarrow \real$ is a given
convex function, and $\Constraint \subset \real^\usedim$ is a convex
constraint set, chosen by the user to enforce a certain type of
structure in the solution.  Important special cases of GLMs include
the linear Gaussian model, in which $\LinkFun(u, y) = \frac{1}{2} (y -
u)^2$, and the problem~\eqref{EqnMLEGLM} corresponds to a regularized
form of least-squares, as well as the problem of logistic regression,
obtained by setting $\LinkFun(u, y)= \log(1+\exp(-yu))$.

Letting $\Amat \in \real^{\numobs \times \usedim}$ denote the data
matrix with $a_i \in \real^\usedim$ as its $i^{th}$ row, the Hessian
of the objective~\eqref{EqnMLEGLM} takes the form
\begin{align*}
\nabla^2 f(x) = \Amat^T \diag\left(\LinkFunDouble(a_i^T
x)\right)_{i=1}^\numobs \Amat 
\end{align*}
Since the function $\LinkFun$ is convex, we are guaranteed that
$\LinkFunDouble(a_i^T x) \geq 0$, and hence the quantity $\diag\left(
\LinkFunDouble(a_i^T x)\right)^{1/2} \Amat$ can be used as an $\numobs
\times \usedim$ matrix square-root.  We return to explore this class
of examples in more depth in Section~\ref{SecConstGLM}.
\end{example}


\subsection{Local convergence analysis using strong convexity}
\label{SecLocal}

Returning now to the general setting, we now begin by proving a local
convergence guarantee for the sketched Newton updates.  In particular,
this theorem provides insight into how large the sketch dimension
$\numproj$ must be in order to guarantee good local behavior of the
sketched Newton algorithm.

This choice of sketch dimension is determined by geometry of the
problem, in particular in terms of the tangent cone defined by the
optimum.  Given a constraint set $\Constraint$ and the minimizer
$\xstar \defn \arg \min \limits_{x \in \Constraint} f(x)$, the tangent
cone at $\xstar$ is given by
\begin{align}
\label{EqnDefnTangentCone}
\ConeSet & \defn \big \{ \Delta \in \real^\usedim \, \mid \, \xstar +
t \Delta \in \Constraint \quad \mbox{for some $t > 0$} \big \}.
\end{align}
Recalling the definition of the Gaussian width from
Section~\ref{SecBackWidths}, our first main result requires the sketch
dimension to satisfy a lower bound of the form
\begin{align}
\label{EqnWidthOne}
\numproj \geq \frac{c}{\epsilon^2} \max_{x\in \Constraint}
\WidthSq(\HessSqrt{x} \ConeSet),
\end{align}
where $\epsilon \in (0,1)$ is a user-defined tolerance, and $c$ is a
universal constant.  Since the Hessian square-root $\HessSqrt{x}$ has
dimensions $\numobs \times \usedim$, this squared Gaussian width is at
at most $\min \{\numobs, \usedim \}$.  This worst-case bound is
achieved for an unconstrained problem (in which case $\ConeSet =
\real^\usedim$), but the Gaussian width can be substantially smaller
for constrained problems.  See the example following
Theorem~\ref{ThmNewtonsSketchStronglyConvex} for an illustration.

In addition to this Gaussian width, our analysis depends on the
cone-constrained eigenvalues of the Hessian $\Hess{\xstar}$, which are
defined as
\begin{align}
\label{EqnConeConstrainedEig}
\relow & = \inf_{z \in \ConeSet \cap \Sphere{\usedim}}
\inprod{z}{\nabla^2 f(\xstar))z}, \quad \mbox{and} \quad \reup =
\sup_{z\in \ConeSet \cap \Sphere{\usedim}} \inprod{z}{\nabla^2
  f(\xstar)) z},
\end{align}
In the unconstrained case ($\Constraint = \real^\usedim$), we have
$\ConeSet = \real^\usedim$, and so that $\relow$ and $\reup$ reduce to
the minimum and maximum eigenvalues of the Hessian $\Hess{\xstar}$.
In the classical analysis of Newton's method, these quantities measure
the strong convexity and smoothness parameters of the function $f$.

With this set-up, the following theorem is applicable to any
twice-differentiable objective $f$ with cone-constrained eigenvalues
$(\relow, \reup)$ defined in equation~\eqref{EqnConeConstrainedEig},
and with Hessian that is $\LipCon$-Lipschitz continuous, as defined in
equation~\eqref{EqnDefnHessianLipCon}.

\btheos[Local convergence of Newton Sketch]
\label{ThmNewtonsSketchStronglyConvex}

For given parameters $\delta, \epsilon \in (0,1)$, consider the Newton
sketch updates~\eqref{EqnNewtonSketch} based on an initialization
$\xit{0}$ such that \mbox{$\|\xit{0} - \xstar\|_2 \leq \delta
  \frac{\relow}{8 \LipCon}$}, and a sketch dimension $\numproj$
satisfying the lower bound~\eqref{EqnWidthOne}.  Then with probability
at least $1-c_1 e^{-c_2 \numproj}$, the $\ell_2$-error satisfies the
recursion
\begin{align}
\label{EqnLocalConvergence}
\| \xit{t+1} - \xstar\|_2 & \leq \epsilon \frac{\reup}{\relow} \|
\xit{t} - \xstar\|_2 + \frac{4 L}{\relow} \|\xit{t} - \xstar\|_2^2.
\end{align}
\etheos

The bound~\eqref{EqnLocalConvergence} shows that when $\epsilon$ is
set to a fixed constant---say $\epsilon = 1/4$---the algorithm
displays a linear-quadratic convergence rate in terms of the error
$\DelIt{t} = \xit{t} - \xstar$.  More specifically, the rate is
initially quadratic---that is, $\|\DelIt{t+1}\|_2 \approx \frac{4
  L}{\relow} \|\DelIt{t}\|_2^2$ when $\|\DelIt{t}\|_2$ is large.
However, as the iterations progress and $\|\DelIt{t}\|_2$ becomes
substantially less than 1, then the rate becomes linear---meaning that
$\|\DelIt{t+1}\|_2 \approx \epsilon \frac{\reup}{\relow}
\|\DelIt{t}\|_2$---since the term $\frac{4L}{\relow}\|\DelIt{t}\|_2^2$
becomes negligible compared to $\epsilon \frac{\reup}{\relow}
\|\DelIt{t}\|_2$.  If we perform $\TotalIt$ steps in total, the linear
rate guarantees the conservative error bounds
\begin{align}
 \| \xit{\TotalIt} - \xstar\|_2 \leq \frac{\relow}{8
   \LipCon}\Big(\frac{1}{2} + \epsilon\frac{\reup}{\relow} \Big)^\TotalIt \, ,
 \quad \mbox{and} \quad f(\xit{\TotalIt}) - f(\xstar) \leq
 \frac{\reup\relow}{8L}\Big(\frac{1}{2} +
 \epsilon\frac{\reup}{\relow}\Big)^\TotalIt\,.
\end{align}

A notable feature of Theorem~\ref{ThmNewtonsSketchStronglyConvex} is
that, depending on the structure of the problem, the linear-quadratic
convergence can be obtained using a sketch dimension $\numproj$ that
is substantially smaller than $\min \{\numobs, \usedim \}$. As an
illustrative example, we performed simulations for some instantiations
of a portfolio optimization problem: it is a linearly-constrained
quadratic program of the form
\begin{align}
\label{EqnPortOne}
\min_{ \substack{x \geq 0 \\\sum_{j=1}^\usedim x_j = 1}} \Big \{
\frac{1}{2} x^T A^T A x - \inprod{c}{x} \Big \},
\end{align}
where $A \in \real^{\numobs \times \usedim}$ and $c \in \real^\usedim$
are empirically estimated matrices and vectors (see
Section~\ref{SecPort} for more details).  We used the Newton sketch to
solve different sizes of this problem \mbox{$\usedim \in \{10, 20, 30,
  40, 50, 60\}$,} and with $\numobs = \usedim^3$ in each case.  Each
problem was constructed so that the optimum $\xstar$ had at most
$\kdim = \lceil 2 \log (\usedim) \rceil$ non-zero entries.  A
calculation of the Gaussian width for this problem (see
Appendix~\ref{AppEllOneWidth} for the details) shows that it suffices
to take a sketch dimension $\numproj \succsim s \log \usedim$, and we
implemented the algorithm with this choice.
\begin{figure}[h]
\begin{center}
\widgraph{.5\textwidth}{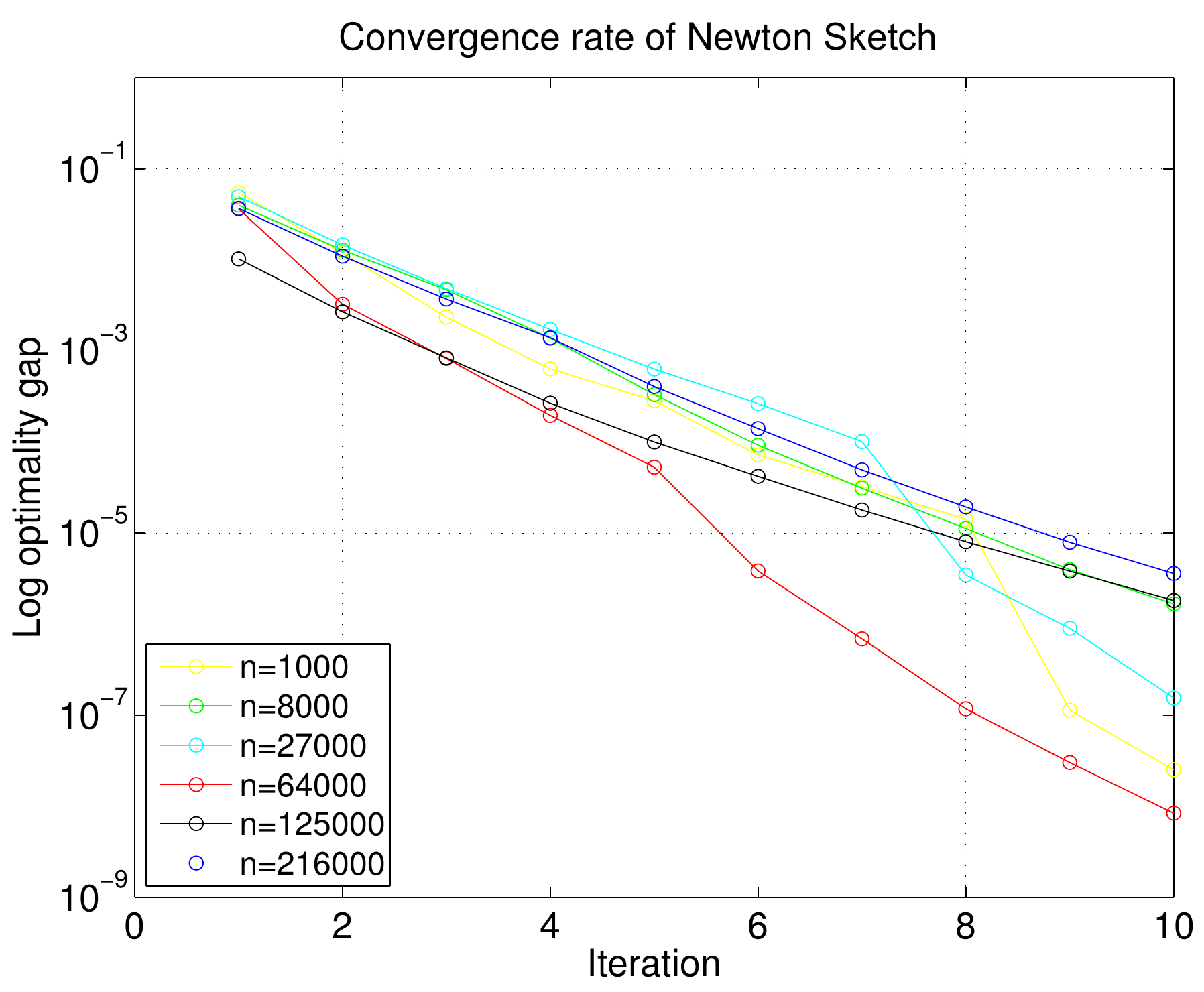}
\end{center}
\caption{Empirical illustration of the linear convergence of the
  Newton sketch algorithm for an ensemble of portfolio optimization
  problems~\eqref{EqnPortOne}.  In all cases, the algorithm was
  implemented using a sketch dimension $\numproj = \lceil 4 s \log
  \usedim \rceil$, where $s$ is an upper bound on the number of
  non-zeros in the optimal solution $x^*$; this quantity satisfies the
  required lower bound~\eqref{EqnWidthOne}, and consistent with the
  theory, the algorithm displays linear convergence.}
\label{FigPort}
\end{figure}
Figure~\ref{FigPort} shows the convergence rate of the Newton sketch
algorithm for the six different problem sizes: consistent with our
theory, the sketch dimension $\numproj \ll \min \{\usedim, \numobs \}$
suffices to guarantee linear convergence in all cases.

It is also possible obtain an asymptotically super-linear rate by
using an iteration-dependent sketching accuracy $\epsilon =
\epsilon(t)$.  The following corollary summarizes one such possible
guarantee:

\bcors 
Consider the Newton sketch iterates using the iteration-dependent
sketching accuracy $\epsilon(t) = \frac{1}{\log(1+t)}$.  Then with the
same probability as in Theorem~\ref{ThmNewtonsSketchStronglyConvex},
we have
\begin{align*}
\| \xit{t+1} - \xstar\|_2 & \leq \frac{1}{\log(1+t)}
\frac{\reup}{\relow} \| \xit{t} - \xstar \|_2 + \frac{4L}{ \relow}
\|\xit{t} - \xstar\|_2^2,
\end{align*}
and consequently, super-linear convergence is obtained---namely,
$\lim_{t\rightarrow\infty}~ \frac{\|\xit{t+1} - \xstar\|_2}{\| \xit{t}
  - \xstar\|_2} = 0$. 
\ecors
\noindent Note that the price for this super-linear convergence is
that the sketch size is inflated by the factor $\epsilon^{-2}(t) =
\log^2(1 + t)$, so it is only logarithmic in the iteration number.


\section{Newton sketch for self-concordant functions}
\label{SecGlobal}

The analysis and complexity estimates given in the previous section
involve the curvature constants $(\relow, \reup)$ and the Lipschitz
constant $\LipCon$, which are seldom known in practice.  Moreover, as
with the analysis of classical Newton method, the theory is local, in
that the linear-quadratic convergence takes place once the iterates
enter a suitable basin of the origin.

In this section, we seek to obtain global convergence results that do
not depend on unknown problem parameters.  As in the classical
analysis, the appropriate setting in which to seek such results is for
self-concordant functions, and using an appropriate form of
backtracking line search.  We begin by analyzing the unconstrained
case, and then discuss extensions to constrained problems with
self-concordant barriers.  In each case, we show that given a suitable
lower bound on the sketch dimension, the sketched Newton updates can
be equipped with global convergence guarantees that hold with
exponentially high probability.  Moreover, the total number of
iterations does not depend on any unknown constants such as strong
convexity and Lipschitz parameters.


\subsection{Unconstrained case}

In this section, we consider the unconstrained optimization problem
$\min_{x \in \real^\usedim} f(x)$, where $f$ is a closed convex
self-concordant function which is bounded below.  Note that a closed
convex function $\phi: \real \rightarrow \real$ is
\emph{self-concordant} if
\begin{align} 
\label{DefnSelfConcordanceSingle} 
|\phi^{\prime \prime \prime}(x)| \le 2 \left(\phi^{\prime
  \prime}(x)\right)^{3/2}.
\end{align}
This definition can be extended to a function $f: \real^\usedim
\rightarrow \real$ by imposing this requirement on the univariate
functions $\phi_{x,y}(t) \defn f(x + t y)$, for all choices of $x,y$
in the domain of $f$. Examples of self-concordant functions include linear and quadratic functions and negative logarithm. Self concordance is preserved under addition and affine transformations.

Our main result provide a bound on the total number of Newton sketch
iterations required to obtain a $\delta$-accurate solution without
imposing any sort of initialization condition (as was done in our
previous analysis).  This bound scales proportionally to
$\log(1/\delta)$ and inversely in a parameter $\NEWGAM$ that depends
on sketching accuracy $\epsilon \in (0,\frac{1}{4})$ and backtracking
parameters $(\linea, \lineb)$ via
\begin{align}
\label{EqnDefnNewGam}
\NEWGAM = \linea \lineb \frac{\eta^2}{1 + (\frac{1 + \epsilon}{1 -
    \epsilon})\eta} \quad \mbox{where} \quad \NEWETA = \frac{1}{8}
\,\frac{1 - \frac{1}{2}(\frac{1 + \epsilon}{1-\epsilon})^2 -
  \linea}{(\frac{1 + \epsilon}{1 - \epsilon})^3}\,.
\end{align}

\begin{algorithm}[h]
\caption{Unconstrained Newton Sketch with backtracking line
  search\label{AlgNewtonSketchBacktracking}}
\begin{algorithmic}[1]
  \scriptsize \REQUIRE Starting point $\xit{0}$, tolerance $\delta >
  0$, $(\linea, \lineb)$ line-search parameters, sketching matrices
  $\{ \SketchIt{t} \}_{t=0}^\infty \in \real^{\numproj \times
    \numobs}$.
\STATE Compute approximate Newton step $\Delta \xit{t}$ and approximate
Newton decrement $\lambda(x)$
\begin{align*}
  \Delta \xit{t} & \defn \arg\min_{\Delta} ~\inprod{\nabla
    f(\xit{t})}{\Delta} + \frac{1}{2} \|S^t(\nabla^2 f(\xit{t}))^{1/2}
  \Delta\|_2^2;\\ \newdecsketch_f(\xit{t}) & \defn \nabla
  f(x)^T\Delta \xit{t}.
\end{align*}
\STATE Quit if $\tilde\newdec(\xit{t})^2/2 \leq \delta$.
\STATE Line search: choose $\mu$ : \quad {\bf while} $f(\xit{t} + \mu
\Delta \xit{t}) > f(\xit{t}) + \linea \mu \lambda(\xit{t}), \quad \mu
\leftarrow \lineb \mu$
\STATE Update: $\xit{t+1} = \xit{t} + \mu \Delta \xit{t}$ 
\ENSURE
minimizer $x^{t}$, optimality gap $\lambda(\xit{t})$
\end{algorithmic}
\end{algorithm}

\btheos
\label{ThmSelfConcordantUnconstrained}
Let $f$ be a strictly convex self-concordant function. Given a
sketching matrix $S\in\real^{m\times n}$ with $m\ge
\frac{c_3}{\epsilon^2} \max_{x\in \Constraint} \rank(\nabla^2 f(x)) =
\frac{c_3}{\epsilon^2}\,d$, the number of total iterations $T$ for
obtaining an $\delta$ approximate solution in function value via
Algorithm~\ref{AlgNewtonSketchBacktracking} is bounded by
\begin{align*}
T = \frac{f(\xit{0})-f(\xstar)}{\NEWGAM} + 0.65
\log_2(\frac{1}{16\delta})\,,
\end{align*}
with probability at least $1 - c_1 \TotalIt e^{-c_2 \numproj}$.
\etheos
The bound in the above theorem shows that the convergence of the
Newton Sketch is independent of the properties of the function $f$ and
problem parameters, similar to classical Newton's method.  Note that
for problems with $\numobs > \usedim$, the complexity of each Newton
sketch step is at most $\order(\usedim^3 + \numobs \usedim \log
\usedim$), which is smaller than that of Newton's Method
($\order(nd^2)$), and also smaller than typical first-order
optimization methods ($\order(\numobs \usedim)$) whenever $\numobs >
\usedim^2$.
%

\subsection{Newton Sketch with self-concordant barriers}

We now turn to the more general constrained case.  Given a closed,
convex self-concordant function $\fbase: \real^\usedim \rightarrow
\real$, let $\Constraint$ be a convex subset of $\real^\usedim$, and
consider the constrained optimization problem $\min_{x \in
  \Constraint} \fbase(x)$.  If we are given a convex self-concordant
barrier function $\gcon$ for the constraint set $\Constraint$, it is
equivalent to consider the unconstrained problem
\begin{align*}
\min_{x \in \real^\usedim} \ENCMIN{\underbrace{f_0(x) + g(x)}_{f(x)}}.
\end{align*}

One way in which to solve this unconstrained problem is by sketching
the Hessian of both $\fzero$ and $\gcon$, in which case the theory of
the previous section is applicable.  However, there are many cases in
which the constraints describing $\Constraint$ are relatively simple,
and so the Hessian of $\gcon$ is highly-structured.  For instance, if
the constraint set is the usual simplex (i.e., $x \geq 0$ and
$\inprod{1}{x} \leq 1$), then the Hessian of the associated log
barrier function is a diagonal matrix plus a rank one matrix.  Other
examples include problems for which $\gcon$ has a separable structure;
such functions frequently arise as regularizers for ill-posed inverse
problems.  Examples of such regularizers include $\ell_2$
regularization $\gcon(x) = \frac{1}{2}\|x\|_2^2$, graph regularization
$\gcon(x)= \frac{1}{2} \sum_{i,j\in {E}} (x_i-x_j)^2$ induced by an
edge set ${E}$ (e.g., finite differences) and also other
differentiable norms $\gcon(x) = \left(\sum_{i=1}^\usedim x_i^p
\right)^{1/p}$ for $1<p<\infty$.

In all such cases, an attractive strategy is to apply a \emph{partial
  Newton sketch}, in which we sketch the Hessian term $\nabla^2
\fzero(x)$ and retain the exact Hessian $\nabla^2 \gcon(x)$, as in the
previously described updates~\eqref{EqnNewtonSketchPartial}.  More
formally, Algorithm~\ref{AlgNewtonSketchBacktrackingBarrier} provides
a summary of the steps, including the choice of the line search
parameters.  The main result of this section provides a guarantee on
this algorithm, assuming that the sequence of sketch dimensions
$\{\numprojit{t}\}_{t=0}^\infty$ is appropriately chosen.

\begin{algorithm}[h]
\caption{Newton Sketch with self-concordant
  barriers \label{AlgNewtonSketchBacktrackingBarrier}}
\begin{algorithmic}[1]
  \scriptsize \REQUIRE Starting point $\xit{0}$, constraint
  $\Constraint$, corresponding barrier function $\gcon$ such that $f =
  \fzero + \gcon$, tolerance $\delta>0$, $(\alpha,\beta)$ line-search
  parameters, sketching matrices $\SketchIt{t} \in \real^{ \numproj
    \times \numobs}$.
\STATE Compute approximate Newton step $\Delta \xit{t}$ and
approximate Newton decrement $\lambda(x)$.
\begin{align*}
  \Delta \xit{t} & \defn \arg\min_{\xit{t}+\Delta \in \Constraint}
  ~\inprod{\nabla f(\xit{t})}{\Delta} + \frac{1}{2} \|S^t(\nabla^2
  \fzero(\xit{t}))^{1/2} \Delta\|_2^2 + \frac{1}{2}\Delta^T \nabla^2
  \gcon(\xit{t}) \Delta ; \\
\newdecsketch_f(\xit{t})& \defn \nabla f(x)^T\Delta \xit{t}
\end{align*}
\STATE Quit if $\tilde \newdec(\xit{t})^2/2\le \delta$.
\STATE Line search: choose $\mu$ : \quad {\bf while} $f(\xit{t} + \mu
\Delta \xit{t}) > f(\xit{t}) + \alpha \mu \lambda(\xit{t}), \quad \mu
\leftarrow \beta\mu$.
\STATE Update: $\xit{t+1} = \xit{t} + \mu \Delta \xit{t}$.  
\ENSURE minimizer $\xit{t}$, optimality gap $\lambda(\xit{t})$.
\end{algorithmic}
\end{algorithm}

The choice of sketch dimensions depends on the tangent cones defined
by the iterates, namely the sets
\begin{align*}
\ConeSetIt{t} \defn \big \{ \Delta \in \real^\usedim \, \mid \,
\xit{t} + \alpha \Delta \in \Constraint \quad \mbox{for some $\alpha >
  0$} \big \}.
\end{align*}
For a given sketch accuracy $\epsilon \in (0,1)$, we require that the
sequence of sketch dimensions satisfies the lower bound
\begin{align}
\label{EqnSketchDimBarrier}
\numprojit{t} \geq \frac{c_3}{\epsilon^2} \max_{x\in \Constraint}
\WidthSq(\nabla^2 f(x)^{1/2} \ConeSetIt{t}).
\end{align}
Finally, the reader should recall the parameter $\NEWGAM$ was defined
in equation~\eqref{EqnDefnNewGam}, which depends only on the sketching
accuracy $\epsilon$ and the line search parameters.  Given this
set-up, we have the following guarantee:

\btheos
\label{ThmSelfConcordantBarrier} 
Let $f:\real^\usedim \rightarrow \real$ be a convex and
self-concordant function, and let \mbox{$g: \real^{\usedim}
  \rightarrow \real \cup \{+\infty\}$} be a convex and self-concordant
barrier for the convex set $\Constraint$.  Suppose that we implement
\mbox{Algorithm~\ref{AlgNewtonSketchBacktrackingBarrier}} with sketch
dimensions $\{\numprojit{t} \}_{t \geq 0}$ satisfying the lower
bound~\eqref{EqnSketchDimBarrier}.  Then taking
\begin{align*}
\TotalIt = \frac{f(\xit{0})-f(\xstar)}{\NEWGAM} + 0.65
\log_2(\frac{1}{16\delta}) \qquad \mbox{iterations},
\end{align*}
suffices to obtain $\delta$-approximate solution in function value
with probability at least $1 - c_1 \TotalIt e^{-c_2 \numproj}$.
\etheos
\noindent Thus, we see that the Newton Sketch method can also be used
with self-concordant barrier functions, which considerably extends its
scope. Section~\ref{SecNumericalLasso} provides a numerical
illustration of its performance in this context. As we discuss in the
next section, there is a flexibility in choosing the decomposition
$f_0$ and $g$ corresponding to objective and barrier, which enables us
to also sketch the constraints.


\subsection{Sketching with interior point methods}

In this section, we discuss the application of Newton Sketch to a form
of barrier or interior point methods.  In particular we discuss two
different strategies and provide rigorous worst-case complexity
results when the functions in the objective and constraints are
self-concordant.  More precisely, let us consider a problem of the
form
\begin{align} 
\label{EqnGeneralConvexConstrainedProblem}
\min_{x \in \real^\usedim} ~   \fzero(x) \quad \mbox{subject
  to} \quad \gcon_j(x) \leq 0 \quad \mbox{for $j= 1,\ldots,
  \CONNUM$,}
\end{align}
where $\fzero$ and $\{\gcon_j\}_{j=1}^\CONNUM$ are
twice-differentiable convex functions.  We assume that there exists a
unique solution $\xstar$ to the above problem.

The barrier method for computing $\xstar$ is based on solving a
sequence of problems of the form
\begin{align} 
\label{EqnBarrierMethod}
\xhat(\tau) \defn \arg \min_{x \in \real^\usedim} \ENCMIN{\tau
  \fzero(x) - \sum_{j=1}^\CONNUM \log(-\gcon_j(x))},
\end{align}
for increasing values of the parameter $\tau \geq 1$.  The family of
solutions $\{\xhat(\tau)\}_{\tau \geq 1}$ trace out what is known as
the central path. A standard bound (e.g.,~\cite{Boyd02}) on the
sub-optimality of $\xhat(\tau)$ is given by
\begin{align*}
\fzero(\xhat(\tau)) - \fzero(\xstar) & \leq \frac{\CONNUM}{\tau}.
\end{align*}
The barrier method successively updates the penalty parameter $\tau$
and also the starting points supplied to Newton's method using
previous solutions. 

Since Newton's method lies at the heart of the barrier method, we can
obtain a fast version by replacing the exact Newton minimization with
the Newton sketch.  Algorithm~\ref{AlgSketchedInteriorPoint} provides
a precise description of this strategy.  As noted in Step 1, there are
two different strategies in dealing with the convex constraints
$\gcon_j(x) \leq 0$ for $j = 1,\ldots,\CONNUM$:
\begin{itemize}
    \item \emph{Full sketch:} Sketch the full Hessian of the objective
      function~\eqref{EqnBarrierMethod} using
      Algorithm~\ref{AlgNewtonSketchBacktracking}\,,
    \item \emph{Partial sketch:} Sketch only the Hessians
      corresponding to a subset of the functions $\{\fzero, \gcon_j, j
      = 1, \ldots, \CONNUM \}$, and use exact Hessians for the other
      functions.  Apply
      Algorithm~\ref{AlgNewtonSketchBacktrackingBarrier}.
\end{itemize}

\begin{algorithm}[t]
\caption{Interior point methods using Newton Sketch}
\label{AlgSketchedInteriorPoint}
\begin{algorithmic}[1]
  \scriptsize \REQUIRE Strictly feasible starting point $\xit{0}$,
  initial parameter $\tau^0$ s.t. $\tau:=\tau^0>0$, $\mu>1$, tolerance
  $\delta>0$.  \STATE Centering step: Compute $\xhat(\tau)$ by Newton
  Sketch with backtracking line-search initialized at $x$ \\using
  Algorithm \ref{AlgNewtonSketchBacktracking} or Algorithm
  \ref{AlgNewtonSketchBacktrackingBarrier}.  \STATE Update
  $x:=\xhat(\tau)$.  \STATE Quit if $r/\tau\le\delta$.  \STATE
  Increase $\tau$ by $\tau:=\mu \tau$.  \ENSURE minimizer
  $\xhat(\tau)$.
\end{algorithmic}
\end{algorithm}

As shown by our theory, either approach leads to the same convergence
guarantees, but the associated computational complexity can vary
depending both on how data enters the objective and constraints, as
well as the Hessian structure arising from particular functions. The
following theorem is an application of the classical results on the
barrier method tailored for Newton Sketch using any of the above
strategies (see e.g.,~\cite{Boyd02}). As before, the key parameter
$\NEWGAM$ was defined in Theorem~\ref{ThmSelfConcordantUnconstrained}.

\btheos[Newton Sketch complexity for interior point methods] 
\label{ThmInterior}
For a given target accuracy $\delta \in (0,1)$ and any $\mu > 1$, the
total number of Newton Sketch iterations required to obtain a
$\delta$-accurate solution using
Algorithm~\ref{AlgSketchedInteriorPoint} is at most
\begin{align}
\ceil[\Bigg]{\frac{\log{(\CONNUM /(\tau^0 \delta)}}{\log \mu
}}\left(\frac{\CONNUM (\mu-1-\log\mu)}{\gamma} + 0.65
\log_2(\frac{1}{16\delta})\right).
\end{align}

\etheos
\noindent 
If the parameter $\mu$ is set to minimize the above upper-bound, the
choice $\mu = 1+\frac{1}{\CONNUM}$ yields $\order(\sqrt{\CONNUM})$
iterations. However, when applying the standard Newton method, this
``optimal'' choice is typically not used in practice: instead, it is
common to use a fixed value of \mbox{$\mu \in [2, 100]$.}  In
experiments, experience suggests that the number of Newton iterations
needed is a constant independent of $\CONNUM$ and other parameters.
Theorem~\ref{ThmInterior} allows us to obtain faster interior point
solvers with rigorous worst-case complexity results.  We show
different applications of Algorithm~\ref{AlgSketchedInteriorPoint} in
the following section.


\section{Applications and numerical results}
\label{SecApplications}

In this section, we discuss some applications of the Newton sketch to
different optimization problems.  In particular, we show various forms
of Hessian structure that arise in applications, and how the Newton
sketch can be computed. When the objective and/or the constraints
contain more than one term, the barrier method with Newton Sketch has
some flexibility in sketching. We discuss the choices of partial
Hessian sketching strategy in the barrier method. It is also possible
to apply the sketch in the primal or dual form, and we provide
illustrations of both strategies here.



\subsection{Estimation in generalized linear models}
\label{SecConstGLM}

Recall the problem of (constrained) maximum likelihood estimation for
a generalized linear model, as previously introduced in
Example~\ref{ExaGLMs}.  It leads to the family of optimization
problems~\eqref{EqnMLEGLM}: here $\LinkFun: \real \rightarrow \real$
is a given convex function arising from the probabilistic model, and
$\Constraint \subseteq \real^\usedim$ is a closed convex set that is
used to enforce a certain type of structure in the solution, Popular
choices of such constraints include $\ell_1$-balls (for enforcing
sparsity in a vector), nuclear norms (for enforcing low-rank structure
in a matrix), and other non-differentiable semi-norms based on total
variation (e.g., $\sum_{j=1}^{\usedim-1} |x_{j+1}-x_j|$), useful for
enforcing smoothness or clustering constraints.

Suppose that we apply the Newton sketch algorithm to the optimization
problem~\eqref{EqnMLEGLM}.  Given the current iterate $\xit{t}$,
computing the next iterate $\xit{t+1}$ requires solving the
constrained quadratic program
\begin{align}
\label{EqnGLMInter}
\min_{x \in \Constraint}~ \left \{ \frac{1}{2}\|\Sketch \diag\left(
\LinkFunDouble( \inprod{a_i}{\xit{t}}, y_i)\right)^{1/2} A (x -
\xit{t}) \|_2^2 + \sum_{i=1}^\numobs
\inprod{x}{\LinkFunSingle(\inprod{a_i}{\xit{t}}, y_i)} \right \} \, .
\end{align}
When the constraint $\Constraint$ is a scaled version of the
$\ell_1$-ball---that is, $\Constraint = \{x \in \real^\usedim \, \mid
\|x\|_1 \leq R \}$ for some radius $R > 0$---the convex
program~\eqref{EqnGLMInter} is an instance of the Lasso
program~\cite{Tibshirani96}, for which there is a very large body of
work.  For small values of $R$, where the cardinality of the solution
$x$ is very small, an effective strategy is to apply a homotopy type
algorithm, also known as LARS~\cite{lars,hastie2007lars}, which solves
the optimality conditions starting from $R=0$. For other sets
$\Constraint$, another popular choice is projected gradient descent,
which is efficient when projection onto $\Constraint$ is
computationally simple.

Focusing on the $\ell_1$-constrained case, let us consider the problem
of choosing a suitable sketch dimension $\numproj$.  Our choice
involves the $\ell_1$-restricted minimal eigenvalue of the data matrix
$\Amat^T \Amat$, which is given by
\begin{align} 
\label{EqnDefnRE}
\sigkminsq(\Amat) \defn \min_{ \substack{\|z\|_2 = 1 \\\|z\|_1\le
    2\sqrt{\kdim}}} \|\Amat z\|^2_2.
\end{align}
Note that we are always guaranteed that $\sigkminsq(\Amat )\ge
\lambda_{\min}(\Amat^T \Amat)$.  It also involves certain quantities
that depend on the function $\LinkFun$, namely
\begin{align*}
\LinkFunDouble_{\min} \defn \min \limits_{x \in \Constraint} \min_{i =
  1, \ldots, \numobs} \, \LinkFunDouble(\inprod{a_i}{x}, y_i), \quad
\mbox{and} \quad \LinkFunDouble_{\max} \defn \max \limits_{x \in
  \Constraint} \max \limits_{i = 1, \ldots, \numobs} \,
\LinkFunDouble( \inprod{a_i}{x},y_i),
\end{align*}
where $a_i \in \real^\usedim$ is the $i^{th}$ row of $\Amat$.  With
this set-up, supposing that the optimal solution $\xstar$ has
cardinality at most $\|\xstar\|_0 \leq \kdim$, then it can be shown
(see Lemma~\ref{LemL1ConstrainedGLMWidth} in
Appendix~\ref{AppEllOneWidth}) that it suffices to take a sketch size
\begin{align}
\label{EqnL1GLMSketchSizeChoice}
\numproj & = c_0 \,
         {\frac{\LinkFunDouble_{\max}}{\LinkFunDouble_{\min}}}
         \frac{\max \limits_{j=1,\ldots, \usedim} \|A_j\|_2^2
         }{\sigkminsq(A)} \; {\kdim \log \usedim},
\end{align}
where $c_0$ is a universal constant.  Let us consider some examples to
illustrate:
\begin{itemize}
  \item Least-Squares regression: $\LinkFun(u) = \frac{1}{2} u^2$,
    $\LinkFunDouble(u)=1$ and $\LinkFunDouble_{\min} =
    \LinkFunDouble_{\max} = 1$.
  \item Poisson regression: $\LinkFun(u) = e^u$, $\LinkFunDouble(u) =
    e^u$ and $\frac{\LinkFunDouble_{\max}}{\LinkFunDouble_{\min}} =
    \frac{e^{R \AMAX}} {e^{-R \AMIN } }$
  \item Logistic regression: $\LinkFun(u) = \log(1+e^u)$,
    $\LinkFunDouble(u) = \frac{e^u}{(e^u+1)^2}$ and
    $\frac{\LinkFunDouble_{\max}}{\LinkFunDouble_{\min}} = \frac{e^{R
      \AMIN}} {e^{-R \AMAX } } \frac{(e^{-R \AMAX } +1)^2}{(e^{R
      \AMIN}+1)^2}$\,,
\end{itemize}
where $\AMAX \defn \max \limits_{i = 1, \ldots, \numobs}
\|a_i\|_\infty$, and $\AMIN \defn \min \limits_{i = 1,\ldots,\numobs}
\|a_i\|_\infty$.

For typical distributions of the data matrices, the sketch size choice
given in equation~\eqref{EqnL1GLMSketchSizeChoice} is $\order(\kdim
\log \usedim)$. As an example, consider data matrices $\Amat \in
\real^{\numobs \times \usedim}$ where each row is independently
sampled from a sub-Gaussian distribution with variance $1$. Then
standard results on random matrices~\cite{Ver11} show that
$\sigkminsq(A) > 1/2$ as long as $\numobs > c_1 \kdim \log \usedim$
for a sufficiently large constant $c_1$.  In addition, we have $\max
\limits_{j = 1, \ldots, \usedim} \|A_j\|_2^2 = \order(\numobs)$, as
well as $\frac{\LinkFunDouble_{\max}}{\LinkFunDouble_{\min}} =
\order(\log(n))$. For such problems, the per iteration complexity of
Newton Sketch update scales as $\order(\kdim^2 \usedim
\log^2(\usedim))$ using standard Lasso solvers
(e.g.,~\cite{kim2007interior}) or as $\order(\kdim \usedim
\log(\usedim))$ using projected gradient descent.  Both of these
scalings are substantially smaller than conventional algorithms that
fail to exploit the small intrinsic dimension of the tangent cone.


\subsection{Semidefinite programs}

The Newton sketch can also be applied to semidefinite programs.  As
one illustration, let us consider the metric learning problem studied
in machine learning. Given feature vectors $a_1,\ldots a_\numobs
\in\real^\usedim$ and corresponding indicator
$y_{ij}\in\{-1,+1\}^\numobs$ where $y_{ij}=+1$ if $a_i$ and $a_j$
belong to the same label and $y_{ij}=-1$ otherwise for all $i\neq j$
and $1\le i,j\le n$. The task is to learn a positive semidefinite
matrix $X$ which represents a metric such that the semi-norm
$\|a\|_{X} \defn \sqrt{\inprod{a}{X a}}$ establishes a nearness
measure depending on class label. Using $\ell_2$-loss, the
optimization can be stated as the following semi-definite program
(SDP)
\begin{align*}
\min_{X\succeq 0} \ENCMIN{ \sum_{i\neq j}^{n\choose 2} ~
  \left(\inprod{X}{(a_i-a_j)(a_i-a_j)^T} - y_{ij}\right)^2 + \lambda
  \trace(X)}.
\end{align*}
Here the term $\trace(X)$, along with its multiplicative pre-factor
$\lambda > 0$ that can be adjusted by the user, is a regularization
term for encouraging a relatively low-rank solution.  Using the
standard self-concordant barrier $X \mapsto \log \det(X)$ for the PSD
cone, the barrier method involves solving a sequence of sub-problems
of the form
\begin{align*}
\min_{X \in \real^{\usedim \times \usedim}} \ENCMIN{\underbrace{\tau
    \sum_{i=1}^\numobs ~ (\inprod{X}{a_i a_i^T} - y_i)^2 + \tau\lambda
    \trace X - \log \det\left(X\right)}_{f(\myvec(X))}}.
\end{align*}
Now the Hessian of the function $\myvec(X) \mapsto f(\myvec(X))$ is a
$\usedim^2 \times \usedim^2$ matrix given by
\begin{align*}
\nabla^2 f \big( \myvec(X) \big) & = \tau \: \sum_{i\neq j}^{\numobs
  \choose 2} \myvec(A_{ij}) \myvec(A_{ij})^T + X^{-1} \otimes X^{-1},
\end{align*}
where $A_{ij} \defn (a_i-a_j)(a_i-a_j)^T$. Then we can apply the
barrier method with partial Hessian sketch on the first term,
$\{S_{ij} \myvec(A_{ij})\}_{i\neq j}$ and exact Hessian for the second
term. Since the vectorized decision variable is $\myvec(X) \in
\real^{\usedim^2}$ the complexity of Newton Sketch is
$\order(\numproj^2 \usedim^2)$ while the complexity of a classical SDP
interior-point solver is $\order(\numobs \usedim^4)$.

\subsection{Portfolio optimization and SVMs}
\label{SecPort}
Here we consider the Markowitz formulation of the portfolio
optimization problem~\cite{Markowitz59}.  The objective is to find $x
\in \real^\usedim$ belonging to the unit simplex, which corresponds to
non-negative weights associated with each of $d$ possible assets, so
as to maximize the expected return minus a coefficient times the
variance of the return.  Letting $\mu \in \real^\usedim$ denote a
vector corresponding to mean return of the assets, and we let $\Sigma
\in \real^{\usedim \times \usedim}$ be a symmetric, positive
semidefinite matrix, covariance of the returns.  The optimization
problem is given by
\begin{align} 
\label{EqnPortfolioOpt}
 \max_{x\ge 0,\,\sum_{j=1}^\usedim x_j \le 1 } \Big \{ \inprod{\mu}{x}
 - \lambda \,x^T \Sigma x \Big \}.
\end{align}
The covariance of returns is often estimated from past stock data via
empirical covariance, $\Sigma = A^T A$ where the columns of $A$ are
time series corresponding to assets normalized by $\sqrt{n}$, where
$n$ is the length of the observation window.

The barrier method can be used solve the above problem by solving
penalized problems of the
\begin{align*}
\min_{x \in \real^\usedim} \Big \{ ~\underbrace{ -\tau\,\mu^T x +
  \tau\lambda \, x^T A^T A x - \sum_{i=1}^{d} \log(\inprod{e_i}{x}) -
  \log(1- \inprod{1}{x})}_{f(x)} \Big \},
\end{align*}
where $e_i \in \real^\usedim$ is the $i^{th}$ element of the canonical
basis and $1$ is row vector of all-ones.  Then the Hessian of the
above barrier penalized formulation can be written as
\begin{align*}
\nabla^2 f(x) = \tau\lambda \, \Amat^T \Amat + \diag
\left(x_i^2\right)^{-1} + 1 1^T
\end{align*}
Consequently we can sketch the data dependent part of the Hessian via
$\tau \lambda \Sketch \Amat$ which has at most rank $\numproj$ and
keep the remaining terms in the Hessian exact. Since the matrix $11^T$
is rank one, the resulting sketched estimate is therefore diagonal
plus rank $(\numproj + 1)$ where the matrix inversion lemma can be
applied for efficient computation of the Newton Sketch update (see
e.g. \cite{Golub96}). Therefore, as long as $\numproj \leq \usedim$,
the complexity per iteration scales as $\order (\numproj \usedim^2)$,
which is cheaper than the $\order(\numobs \usedim^2)$ per step
complexity associated with classical interior point methods. We also
note that support vector machine classification problems with squared
hinge loss also has the same form as in \eqref{EqnPortfolioOpt} (see
e.g. \cite{PilWai14a}) where the same strategy can be applied.


\subsection{Unconstrained logistic regression with $\usedim \ll \numobs$}

\begin{figure}[h!]
\begin{center}
\widgraph{0.9\textwidth}{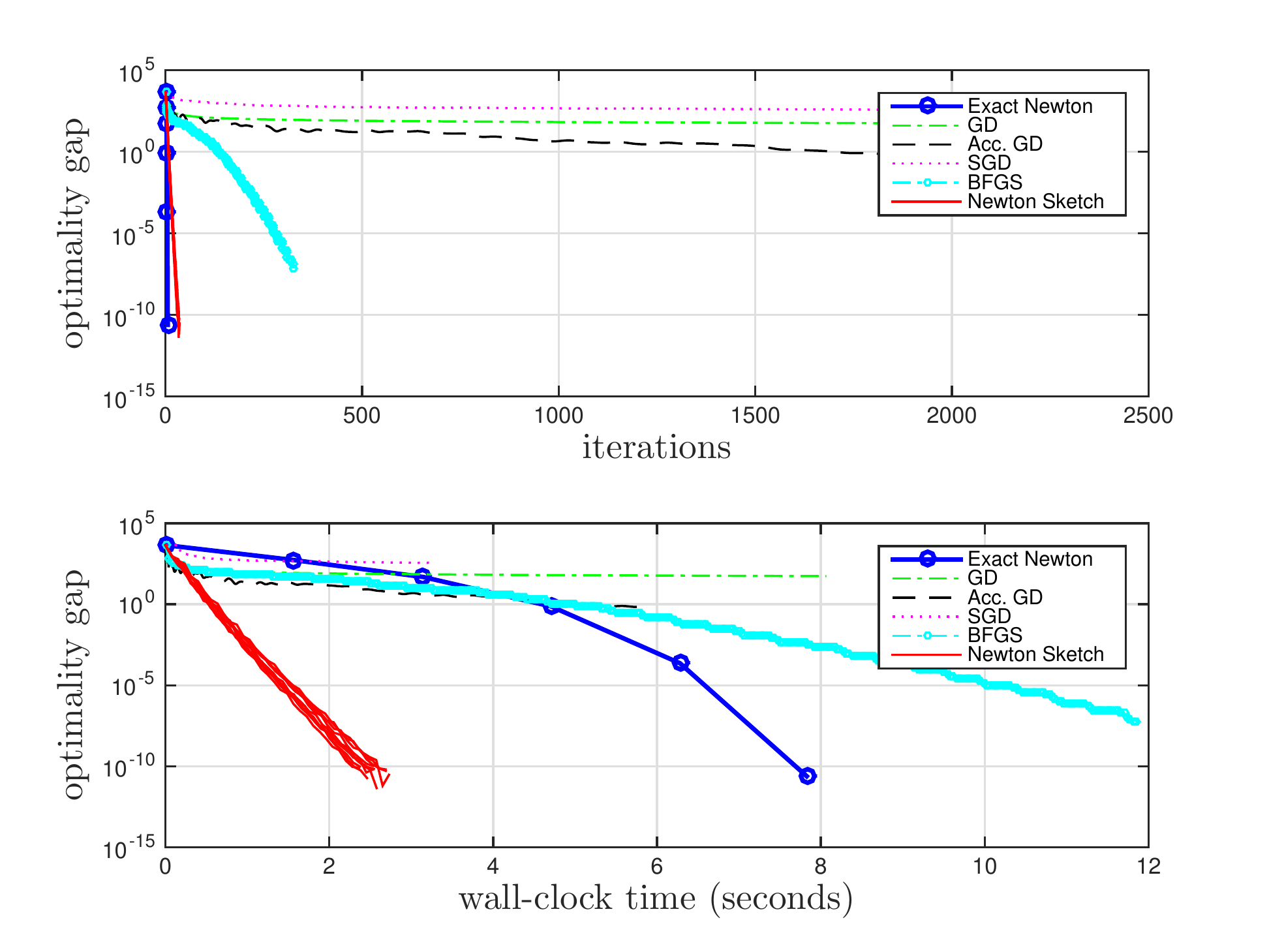}
\caption{Newton Sketch algorithm outperforms other popular
  optimization methods. Plots of the log optimality gap versus
  iteration number (top) and plots of the log optimality gap versus
  wall-clock time (bottom). Newton Sketch empirically provides the
  best accuracy in smallest wall-clock time, and does not require
  knowledge of problem-dependent quantities (such as strong convexity
  and smoothness parameters).}
\label{FigLogistic}
\end{center}
\end{figure}

Let us now turn to some numerical comparisons of the Newton Sketch
with other popular optimization methods for large-scale instances of
logistic regression. More specifically, we generated a feature matrix
$\Amat \in \real^{\numobs \times \usedim}$ based on $\usedim = 100$
features and $\numobs = 16384$ observations. Each row $a_i \in
\real^{\usedim}$ was generated from the $\usedim$-variate Gaussian
distribution $N(0,\Sigma)$ where $\Sigma_{ij}=2 |0.99|^{i-j}$. As
shown in Figure~\ref{FigLogistic}, the convergence of the algorithm
per iteration is very similar to Newton's method. Besides the original
Newton's method, the other algorithms compared are
\begin{itemize}
  \item Gradient Descent (GD) with backtracking line search
  \item Accelerated Gradient Descent (Acc. GD) adapted for strongly
    convex functions with manually tuned parameters.
  \item Stochastic Gradient Descent (SGD) with the classical step size
    choice $1/\sqrt{t}$
  \item Broyden–-Fletcher–-Goldfarb-–Shanno algorithm (BFGS)
    approximating the Hessian with gradients.
\end{itemize}
For each problem, we averaged the performance of the randomized
algorithms (Newton sketch and SGD) over $10$ independent trials.  We
ran the Newton sketch algorithm with sketch size $\numproj = 6
\usedim$.  To be fair in comparisons, we performed hand-tuning of the
stepsize parameters in the gradient-based methods so as to optimize
their performance.  The top panel in Figure~\ref{FigLogistic} plots
the log duality gap versus the number of iterations: as expected, on
this scale, the classical form of Newton's method is the fastest,
whereas the SGD method is the slowest.  However, when the log
optimality gap is plotted versus the wall-clock time in the bottom
panel, we now see that the Newton sketch is the fastest.

\subsection{A dual example: Lasso with $\usedim \gg \numobs$}
\label{SecNumericalLasso}

The regularized Lasso problem takes the form $\min \limits_{x \in
  \real^\usedim} \big \{ \frac{1}{2} \,\|Ax-y\|_2^2 + \lambda \|x\|_1
\big \}$, where $\lambda > 0$ is a user-specified regularization
parameter.  In this section, we consider efficient sketching
strategies for this class of problems in the regime $\usedim \gg
\numobs$.  In particular, let us consider the corresponding dual
program, given by
\begin{align*}
\max_{\|A^T w\|_{\infty} \le \lambda } \ENCMIN{ -\frac{1}{2} \,
  \|y-w\|_2^2}.
\end{align*}
By construction, the number of constraints $\usedim$ in the dual
program is larger than the number of optimization variables $\numobs$.
If we apply the barrier method to solve this dual formulation, then we
need to solve a sequence of problems of the form
\begin{align*}
\min_{w \in \real^\numobs} \ENCMIN{\underbrace{\tau \|y-w\|_2^2 -
    \sum_{j=1}^\usedim \log(\lambda - \inprod{A_j}{w}) -
    \sum_{j=1}^\usedim \log(\lambda + \inprod{A_j}{w})}_{f(x)}},
\end{align*}
where $A_j \in \real^\numobs$ denotes the $j^{th}$ column of $\Amat$.
The Hessian of the above barrier penalized formulation can be written
as
\begin{align*}
\nabla^2 f(w) = \tau I_\numobs + \Amat \diag \left(\frac{1}{(\lambda -
  \inprod{A_j}{w})^2}\right) \Amat^T + \Amat \diag
\left(\frac{1}{(\lambda + \inprod{A_j}{w})^2}\right) \Amat^T\,,
\end{align*}
Consequently we can keep the first term in the Hessian, $\tau I$ exact
and apply partial sketching to the Hessians of the last two terms via
\begin{align*}
\Sketch \diag \left( \frac{1}{|\lambda - \inprod{A_j}{w}|} +
\frac{1}{|\lambda + \inprod{A_j}{w}|} \right) \Amat^T\,.
\end{align*}
Since the partially sketched Hessian is of the form $t I_\numobs + V
V^T$, where $V$ is rank at most $\numproj$, we can use matrix
inversion lemma for efficiently calculating Newton Sketch updates. The
complexity of the above strategy for $\usedim > \numobs$ is
$\order(\usedim \numproj^2)$, where $\numproj$ is at most $\usedim$,
whereas traditional interior point solvers are typically
$\order(\usedim \numobs^2)$ per iteration.

\begin{figure}[h!]
\begin{center}
\begin{tabular}{c}
\widgraph{.8\textwidth}{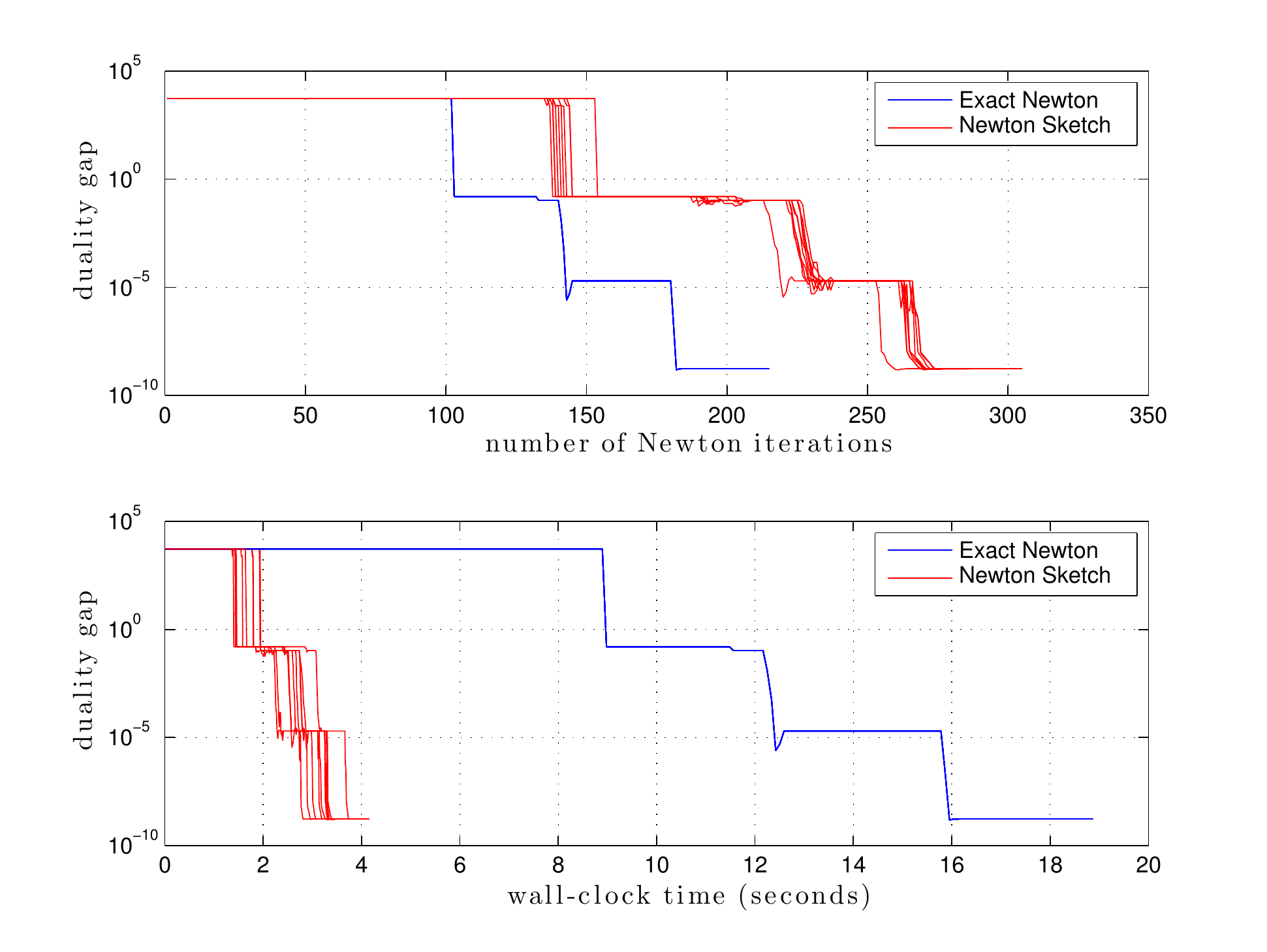} 
\end{tabular}
\caption{Plots of the duality gap versus iteration number (top panel)
  and duality gap versus wall-clock time (bottom panel) for the
  original barrier method (blue) and sketched barrier method (red).
  The sketched interior point method is run $10$ times independently
  yielding slightly different curves in red.  While the sketched
  method requires more iterations, its overall wall-clock time is much
  smaller.}
\label{FigLassoDualityGap}
\end{center}
\end{figure}

\begin{figure}[h!]
\begin{center}
\begin{tabular}{c}
\widgraph{.9\textwidth}{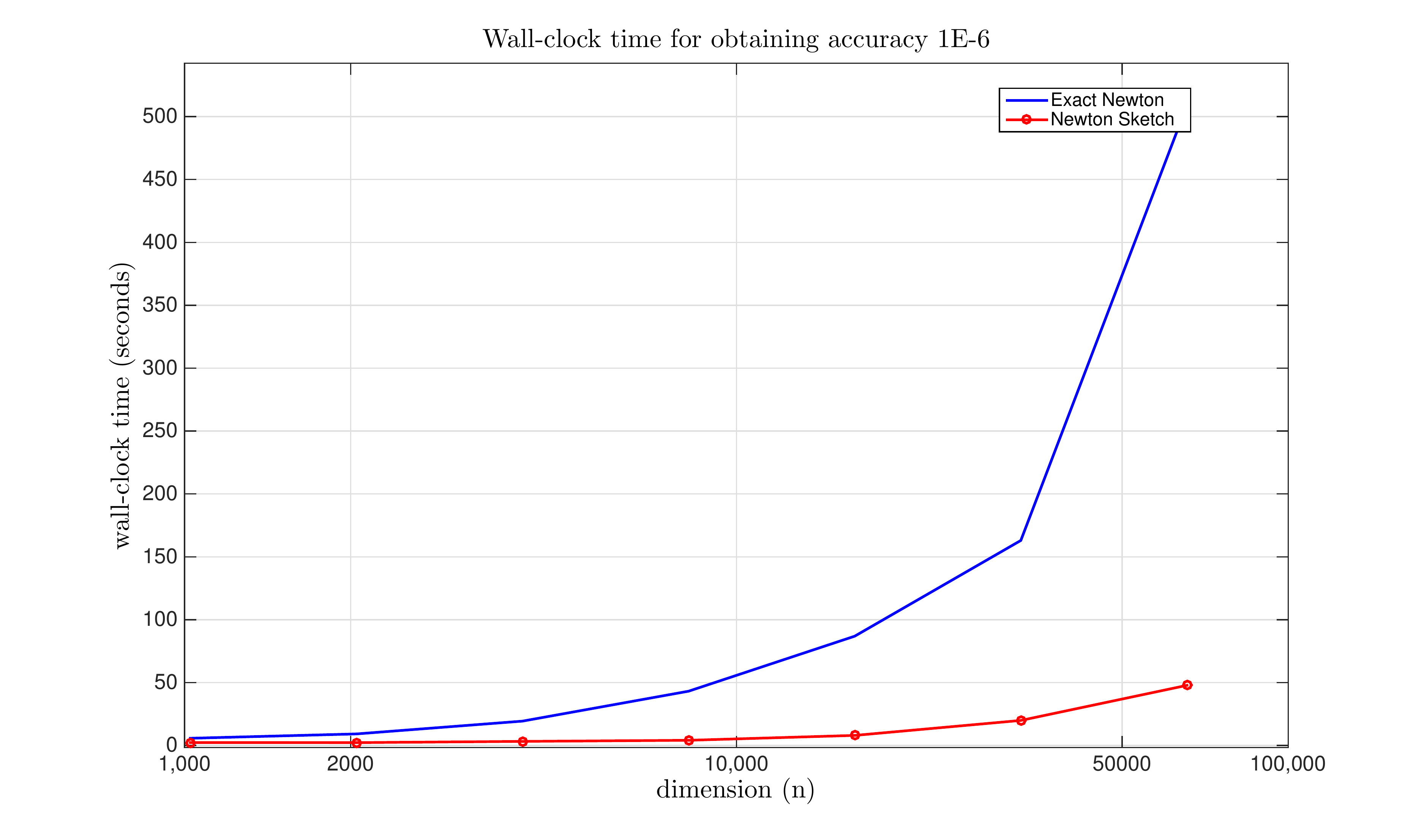} \\
\end{tabular}
\caption{Plot of the wall-clock time in seconds for reaching a duality
  gap of $10^{-6}$ for the standard and sketched interior point
  methods as $\numobs$ increases (in log-scale). The sketched interior
  point method has significantly lower computation time compared to
  the original method.}
\label{FigLassoWallClock}
\end{center}
\end{figure}

In order to test this algorithm, we generated a feature matrix $\Amat
\in \real^{\numobs \times \usedim}$ with $\usedim = 4096$ features and
$\numobs = 50$ observations.  Each row $a_i \in \real^{\usedim}$ was
generated from the multivariate Gaussian distribution $N(0,\Sigma)$
with $\Sigma_{ij}=2* |0.99|^{i-j}$.  For a given problem instance, we
ran $10$ independent trials of the sketched barrier method, and
compared the results to the original barrier method.
Figure~\ref{FigLassoDualityGap} plots the the duality gap versus
iteration number (top panel) and versus the wall-clock time (bottom
panel) for the original barrier method (blue) and sketched barrier
method (red): although the sketched algorithm requires more
iterations, these iterations are cheaper, leading to a smaller
wall-clock time.  This point is reinforced by
Figure~\ref{FigLassoWallClock}, where we plot the wall-clock time
required to reach a duality gap of $10^{-6}$ versus the number of
features $\numobs$ in problem families of increasing size.  Note that
the sketched barrier method outperforms the original barrier method,
with significantly less computation time for obtaining similar
accuracy.

\section{Proofs}
\label{SecProofs}

We now turn to the proofs of our theorems, with more technical details
deferred to the appendices.


\subsection{Proof of Theorem~\ref{ThmNewtonsSketchStronglyConvex}}

Throughout this proof, we let $r \in \Sphere{\usedim}$ denote a fixed
vector that is independent of the sketch matrix $\SketchIt{t}$ and the
current iterate $\xit{t}$.  We then define the following pair of
random variables
\begin{align*}
\ZONEHACK{\Sketch}{x} & \defn \sup_{w \in \HackCone{x} \cap
  \Sphere{\numobs}} \inprod{w}{ \big(\Sketch^T \Sketch -I \big) r}, \\
\ZTWOHACK{\Sketch}{x} & \defn \inf_{w \in \HackCone{x} \cap
  \Sphere{\numobs}} \|\Sketch w\|_2^2.
\end{align*}
These random variables are significant, because the core of our proof
is based on establishing that the error vector $\DelIt{t} = \xit{t} -
\xstar$ satisfies the recursive bound
\begin{align}
\label{EqnKeyRecursionOne}
 \| \DelIt{t+1}\|_2 & \leq \frac{6 \reup \, \ZONET }{\relow \; \ZTWOT}
 \| \DelIt{t}\|_2 + \frac{4 \LipCon}{\relow \; \ZTWOT}
 \|\DelIt{t}\|_2^2,
\end{align}
where $\ZONET \defn \ZONEHACK{\SketchIt{t}}{\xit{t}}$ and $\ZTWOT
\defn \ZTWOHACK{\SketchIt{t}}{\xit{t}}$.
 We then combine this recursion with the
following probabilistic guarantee on $\ZONET$ and $\ZTWOT$.
For a given tolerance parameter $\epsilon \in (0,
\frac{1}{2}]$, consider the ''good event''
\begin{align}
\Event^t & \defn \biggr \{ \ZINF^t \leq \frac{\epsilon}{2}, \mbox{ and }
\ZSUP^t \geq 1-\epsilon \biggr \}.
\end{align}

\blems[Sufficient conditions on sketch dimension~\cite{PilWai14a}]
\label{LemKeyProbabilistic}
\hfill
\begin{enumerate}
\item[(a)] For sub-Gaussian sketch matrices, given a sketch size
  $\numproj > \frac{c_0}{\epsilon^2} \max_{x\in\Constraint} \Width^2(\Hess{x}^{1/2}\ConeSet)$, we have
\begin{align}
\label{EqnGaussianSketch}
\mprob \big[ \Event^t] & \geq 1 - c_1 \CEXP{-c_2 \numproj
  \epsilon^2}.
\end{align}
\item[(b)] For randomized orthogonal system (ROS) sketches over the
  class of self-bounding cones, given a sketch size $\numproj >
  \frac{c_0 \, \log^4 n }{\epsilon^2} \max_{x \in \Constraint} \Width^2(\Hess{x}^{1/2}\ConeSet)$, we
  have
\begin{align}
\mprob \big[ \Event^t] & \geq 1 - c_1 \CEXP{-c_2 \frac{\numproj
    \epsilon^2}{\log^4 n }}.
\end{align}
\end{enumerate}
\elems

\noindent Combining Lemma~\ref{LemKeyProbabilistic} with the
recursion~\eqref{EqnKeyRecursionOne} and re-scaling $\epsilon$ appropriately yields the claim of the theorem.
Accordingly, it remains to prove the
recursion~\eqref{EqnKeyRecursionOne}, and we do so via a basic
inequality argument.  Recall the function $x \mapsto
\SPECFUNTWO{x}{\SketchIt{t}}$ that underlies the sketch Newton
update~\eqref{EqnNewtonSketch}: since $\xit{t}$ and $\xstar$ are
optimal and feasible for the constrained optimization problem, we have
$\SPECFUNTWO{x}{\SketchIt{t}} \leq \SPECFUNTWO{\xstar}{\SketchIt{t}}$.
Introducing the error vector $\DelIt{t} \defn \xit{t} - \xstar$, some
straightforward algebra then then leads to the \emph{basic inequality}
\begin{align}
\label{EqnBasicStart}
\frac{1}{2} \|\SketchIt{t} \HessSqrt{\xit{t}} \DelIt{t+1}\|_2^2 & \leq
\inprod{\SketchIt{t} \hessh \DelIt{t+1}}{S\hessh \DelIt{t}} -
\inprod{\nabla f(\xit{t})-\nabla f(\xstar)}{\DelIt{t+1}}
\end{align}
Let us first upper bound the right-hand side.  By using the integral
form of Taylor's expansion, we have $\inprod{\nabla f(\xit{t})-\nabla
  f(\xstar)}{\DelIt{t+1}} = \int_0^1 \inprod{\Hess{\xit{t} + u(\xstar
    - \xit{t})} \DelIt{t}} {\DelIt{t+1}}du$, and hence
\begin{align*}
\RHS & = \int_0^1 \inprod{ \left[ \hessh (\SketchIt{t})^T \SketchIt{t}
    \hessh - \Hess{\xit{t} + u(\xstar - \xit{t})} \right]\DelIt{t}}
     {\DelIt{t+1}}du
\end{align*}
By adding and subtracting terms and then applying triangle inequality,
we have the bound $\RHS \leq T_1 + T_2$, where
\begin{align*}
\Term_1 & \defn \Big| \int_0^1 \inprod{ \left[ \hessh \big(
    (\SketchIt{t})^T \SketchIt{t} - I \big) \hessh \right]
  \DelIt{t}}{\DelIt{t+1}} \Big|, \quad \mbox{and} \\
\Term_2 & \defn \int_0^1 \opnorm{\Hess{\xit{t} + u (\xstar - \xit{t})}
  - \Hess{\xit{t}}} du \; \|\DelIt{t}\|_2 \|\DelIt{t+1}\|_2.
\end{align*}
Now observe that the vector $r \defn \hessh \DelIt{t}$ is independent
of the randomness in $\SketchIt{t}$, whereas the vector $\hessh
\DelIt{t+1}$ belongs to the cone $\HackCone{\xit{t}}$.  Consequently,
by the definition of $\Zone$, we have
\begin{align}
\label{EqnHessianSmoothnessBound}
\Term_1 &\leq \ZONET \|\hessh \DelIt{t}\|_2 \|\hessh \DelIt{t+1}\|_2\,.
\end{align}
Now note that using the fact that $\reup$ controls the smoothness of
the gradient and the Lipschitz continuity of Hessian we can upper
bound the terms on the above right-hand side as follows
\begin{align*}
\inprod{\DelIt{t}}{\Hess{\xstar} \DelIt{t}} & = \inprod{\DelIt{t}}{
  \Hess{\xstar} \DelIt{t}} + \inprod{\DelIt{t}}{
  \left(\Hess{\xit{t}}-\Hess{\xstar}\right) \DelIt{t}} \\
& \leq \left\{ \beta + \LipCon\|\DelIt{t}\|_2 \right\}
\|\DelIt{t}\|_2^2\,,
\end{align*}
and similarly, $\inprod{\DelIt{t+1}}{\Hess{\xstar} \DelIt{t+1}} \le
\left \{ \beta + \LipCon \|\DelIt{t}\|_2\right\} \|\DelIt{t+1}\|_2^2$.
Combining the above bounds with \eqref{EqnHessianSmoothnessBound} we
obtain
\begin{align}
\Term_1 &\leq \ZONET \left \{ \reup + \LipCon \|\DelIt{t}\|_2\right\} \|\DelIt{t+1}\|_2 \|\DelIt{t}\|_2\,.
\end{align}
On the other hand, by the $\LipCon$-Lipschitz condition on the
Hessian, we have
\begin{align*}
\Term_2 \leq \LipCon \|\DelIt{t}\|_2^2 \, \|\DelIt{t+1}\|_2.
\end{align*}
Substituting these two bounds into our basic inequality, we have
\begin{align}
\label{EqnBasicIntermediate}
\frac{1}{2} \|\SketchIt{t} \HessSqrt{\xit{t}} \DelIt{t+1}\|_2^2 & \leq
\ZONET \,  \left \{ \reup + \LipCon \|\DelIt{t}\|_2\right\} \|\DelIt{t}\|_2 \, \|\DelIt{t+1}\|_2 + \LipCon
\|\DelIt{t}\|_2^2 \, \|\DelIt{t+1}\|_2.
\end{align}
Our final step is to lower bound the left-hand side (LHS) of this
inequality.  By definition of $\Ztwo$, we have
\begin{align*}
\frac{\|\SketchIt{t} \HessSqrt{\xit{t}} \DelIt{t+1}\|^2_2}{ \ZTWOT} &
\geq \inprod{\DelIt{t+1}}{\Hess{\xit{t}} \DelIt{t+1}} \\
& = \inprod{\DelIt{t+1}}{\Hess{\xstar} \DelIt{t+1}} +
\inprod{\DelIt{t+1}}{\big(\Hess{\xit{t}} - \Hess{\xstar}
  \big)\DelIt{t+1}}\\
& \geq \Big \{ \relow - \LipCon \|\DelIt{t}\|_2 \Big \}
\|\DelIt{t+1}\|^2_2.
\end{align*}
Substituting this lower bound into the previous
inequality~\eqref{EqnBasicIntermediate} and then rearranging, we find
that, as long as $\|\DelIt{t}\|_2 < \frac{\relow}{2 \LipCon}$, we also have $\|\DelIt{t}\|_2 < \frac{\reup}{2 \LipCon}$ and consequently
\begin{align*}
\|\DelIt{t+1}\|_2 & \leq \frac{2 \ZONET  \big( \reup + \LipCon \|\DelIt{t}\|_2\big) \, }{ \ZTWOT \big(\relow -
  \LipCon \|\DelIt{t}\|_2 \big) } \|\DelIt{t}\|_2 \, + \frac{2
  \LipCon}{ \big(\relow - \LipCon \|\DelIt{t}\|_2 \big) \, \ZTWOT}
\|\DelIt{t}\|_2^2 \\
& \leq \frac{6 \reup \, \ZONET \, }{\relow \; \ZTWOT} \|
\DelIt{t}\|_2 + \frac{4 \LipCon}{\relow \; \ZTWOT} \|\DelIt{t}\|_2^2,
\end{align*}
as claimed.


\subsection{Proof of Theorem~\ref{ThmSelfConcordantUnconstrained}}

Recall that in this case, we assume that $f$ is a self-concordant
strictly convex function. We adopt the following notation and
conventions from the book~\cite{NesNem94}.  For a given $x \in
\real^\usedim$, we define the pair of dual norms
\begin{align*}
\|u\|_x \, \defn \, \inprod{\fpp(x)u}{u}^{1/2}, \quad \mbox{and} \quad
\|v\|_x^* \, \defn \, \inprod{\fpp(x)^{-1}v}{v}^{1/2}\,,
\end{align*}
as well as the Newton decrement
\begin{align*}
\newdec_f(x) = \inprod{\fpp(x)^{-1}\fp(x)}{\fp(x)}^{1/2} &=
\|\fpp(x)^{-1}\fp(x)\|_x \; = \; \|\fpp(x)^{-1/2}\fp(x)\|_2\,.
\end{align*}
Note that $\fpp(x)^{-1}$ is well-defined for strictly convex
self-concordant functions.  In terms of this notation, the exact
Newton update is given by $x \mapsto \xnewt \defn x + \vnew$, where
\begin{align}
\label{EqnNewtonOrig}
\vnewt & \defn \arg \min_{z \in \Constraint-x}
\ENCMIN{\underbrace{\frac{1}{2}\| \fpp(x)^{1/2} z\|_2^2 +
    \inprod{z}{\nabla f(x)}}_{\Phi(z)}},
\end{align}
whereas the Newton sketch update is given by $x \mapsto \xsketch \defn
x + \vsketch$, where
\begin{align}
\label{EqnNewtonApprox}
\vsketch & \defn \arg \min_{z \in \Constraint-x} \ENCMIN{
  \frac{1}{2}\| S \fpp(x)^{1/2} z\|_2^2 + \inprod{z}{\nabla f(x)}} \,.
\end{align}
The proof of Theorem~\ref{ThmSelfConcordantUnconstrained} given in
this section involves the unconstrained case ($\Constraint =
\real^\usedim$), whereas the proofs of later theorems involve the more
general constrained case.  In the unconstrained case, the two updates
take the simpler forms
\begin{align*}
\xnewt = x - (\Hess{x})^{-1} \nabla f(x), \qquad \mbox{and} \qquad
\xsketch = x - (\HessSqrt{x} \Sketch^T \Sketch \HessSqrt{x})^{-1}
\nabla f(x).
\end{align*}

For a self-concordant function, the sub-optimality of the Newton
iterate $\xnewt$ in function value satisfies the bound
\begin{align*}
f(\xnewt) - \underbrace{\min_{x \in \real^\usedim} f(x)}_{f(\xstar)}
\leq \big[\newdec_f(\xnewt) \big]^2~.
\end{align*}
This classical bound is not directly applicable to the Newton sketch
update, since it involves the \emph{approximate} Newton decrement
\mbox{$\newdecsketch_f(x) \defn -\inprod{\nabla f(x)}{\vsketch}$,} as
opposed to the exact one $\newdec_f(x) \defn -\inprod{\nabla
  f(x)}{\vnewt}$.  Thus, our strategy is to prove that with high
probability over the randomness in the sketch matrix, the approximate
Newton decrement can be used as an exit condition.

Recall the definitions~\eqref{EqnNewtonOrig}
and~\eqref{EqnNewtonApprox} of the exact $\vnewt$ and sketched Newton
$\vsketch$ update directions, as well as the definition of the tangent
cone $\KCONE$ at $x \in \Constraint$. Let $\KCONE^t$ be the tangent
cone at $x^t$. The following lemma provides a high probability bound
on their difference:
\blems 
\label{LemApproxNewtonStep} 
Let $\Sketch \in \real^{m\times n}$ be a sub-Gaussian or ROS sketch
matrix, and consider any fixed vector $x \in \Constraint$ independent
of the sketch matrix.  If $\numproj \ge c_0\frac{\Width(\fpp(x)^{1/2}
  \KCONE^t)^2}{\epsilon^2}$, then
\begin{align} 
\label{EqnLemSketch}
\left\| \fpp(x)^{1/2} (\vsketch-\vnewt) \right\|_2 \leq \epsilon
\left\| \fpp(x)^{1/2} \vnewt \right\|_2~
\end{align}
with probability at least $1-c_1e^{-c_2 \numproj \epsilon^2}$.
\elems

Similar to the standard analysis of Newton's method, our analysis of
the Newton sketch algorithm is split into two phases defined by the
magnitude of the decrement $\newdecsketch_f(x)$.  In particular, the
following lemma constitute the core of our proof:
\blems
\label{LemKeyDecrement}
For $\epsilon \in (0, 1/2)$, there exist constants $\NEWGAM > 0$ and
$\NEWETA \in (0, 1/16)$ such that:
\begin{enumerate}
\item[(a)] If $\newdecsketch_f(x) > \NEWETA$, then $f(\xsketch)-f(x)
  \leq - \NEWGAM$ with probability at least $1-c_1 e^{-c_2 \numproj
    \epsilon^2}$.
\item[(b)] Conversely, if $\newdecsketch_f(x) \leq \NEWETA$, then 
\begin{subequations}
\begin{align}
\label{EqnSecondPhaseControl}
\newdecsketch_f(\xsketch) & \leq \newdecsketch_f(x), \quad \mbox{and} \\
\label{EqnSecondPhaseContract}
\newdec_f(\xsketch) & \leq \big(\frac{16}{25} \big) \newdec_f(x),
\end{align}
\end{subequations}
where both bounds hold with probability $1 - c_1 e^{c_2 \numproj
  \epsilon^2}$.
\end{enumerate}
\elems

\noindent Using this lemma, let us now complete the proof of the
theorem, dividing our analysis into the two phases of the algorithm.

\paragraph{First phase analysis:}
By Lemma~\ref{LemKeyDecrement}(a) each iteration in the first phase
decreases the function value by at least $\NEWGAM > 0$, the number of
first phase iterations $\TotalIt_1$ is at most
\begin{align*}
\TotalIt_1 \defn \frac{f(\xit{0})- f(\xstar)}{\NEWGAM}\,,
\end{align*}
with probability at least $1- \TotalIt_1 c_1 e^{-c_2 \numproj}$.

\paragraph{Second phase analysis:}

Next, let us suppose that at some iteration $t$, the condition
$\newdecsketch_f(\xit{t}) \leq \NEWETA$ holds, so that part (b) of
Lemma~\ref{LemKeyDecrement} can be applied.  In fact, the
bound~\eqref{EqnSecondPhaseControl} then guarantees that
$\newdecsketch_f(\xit{t+1}) \leq \NEWETA$, so that we may apply the
contraction bound~\eqref{EqnSecondPhaseContract} repeatedly for
$\TotalIt_2$ rounds so as to obtain that
\begin{align*}
\newdec_f(\xit{t+\TotalIt_2}) & \leq \big(\frac{16}{25}
\big)^{\TotalIt_2} \newdec_f(\xit{t})
\end{align*}
with probability $1 - \TotalIt_2 c_1 e^{c_2 \numproj}$.

Since $\newdec_f(\xit{t}) \leq \NEWETA \leq 1/16$ by assumption, the
self-concordance of $f$ then implies that
\begin{align*}
f(\xit{t+k}) - f(\xstar) \leq \left(\frac{16}{25}\right)^{k}
\frac{1}{16}.
\end{align*}
Therefore, in order to ensure that and consequently for achieving
$f(\xit{t + k})-f(\xstar) \leq \epsilon$, it suffices to the number of
second phase iterations lower bounded as $\TotalIt_2 \geq 0.65
\log_2(\frac{1}{16\epsilon})$.

Putting together the two phases, we conclude that the total number of
iterations $\TotalIt$ required to achieve $\epsilon$- accuracy is
at most
\begin{align*}
\TotalIt & = \TotalIt_1 + \TotalIt_2 \; \leq \;
\frac{f(\xit{0})-f(\xstar)}{\gamma} + 0.65
\log_2(\frac{1}{16\epsilon})\,,
\end{align*}
and moreover, this guarantee holds with probability at least $1-
\TotalIt c_1 e^{-c_2 \numproj \epsilon^2}$.


\vspace*{.25in}

The final step in our proof of the theorem is to establish
Lemma~\ref{LemKeyDecrement}, and we do in the next two subsections.


\subsubsection{Proof of Lemma~\ref{LemKeyDecrement}(a)}

Our proof of this part is performed conditionally on the event
$\Devent \defn \{\newdecsketch_f(x) > \NEWETA\}$.  Our strategy is to
show that the backtracking line search leads to a stepsize $\stepsize
> 0$ such that function decrement in moving from the current iterate
$x$ to the new sketched iterate $\xsketch = x + \stepsize \vsketch$ is
at least
\begin{align}
\label{EqnDampedPhase}
f(\xsketch) - f(x) \leq - \NEWGAM \quad \mbox{with probability at
  least $1-c_1 e^{-c_2 \numproj}$.}
\end{align}

The outline of our proof is as follows.  Defining the univariate
function $\gfun(u) \defn f(x + u \vsketch)$ and $\epsilon' = \frac{2
  \epsilon}{1 - \epsilon}$, we first show that $\shat = \frac{1}{1 +
  (1 + \epsilon') \newdecsketch_f(x)}$ satisfies the bound
\begin{subequations}
\begin{align}
\label{EqnDampedDecrementExit}
\gfun(\shat) & \leq \gfun(0) - \linea \shat \newdecsketch_f(x)^2,
\end{align}
which implies that $\shat$ satisfies the exit condition of
backtracking line search. Therefore, the stepsize $\stepsize$ must be
lower bounded as $\stepsize \geq \lineb \shat$, which then implies that
the updated solution $\xsketch = x + \stepsize \vsketch$ satisfies the
decrement bound
\begin{align}
\label{EqnDampedDecrementExitTwo}
f(\xsketch) - f(x) & \leq - \linea \lineb \frac{
  \newdecsketch_f(x)^2}{1+(1 + \frac{2 \epsilon}{1-\epsilon})
  \newdecsketch_f(x)}.
\end{align}
\end{subequations}
Since $\newdecsketch_f(x) > \NEWETA$ by assumption and the function $u
\rightarrow\frac{u^2}{1+(1+ \frac{2 \epsilon}{1-\epsilon}) u}$ is
monotone increasing, this bound implies that
inequality~\eqref{EqnDampedPhase} holds with \mbox{$\NEWGAM = \linea
  \lineb \frac{ \NEWETA^2}{1+(1 + \frac{2 \epsilon}{1-\epsilon})
    \NEWETA}$.}\\

It remains to prove the claims~\eqref{EqnDampedDecrementExit}
and~\eqref{EqnDampedDecrementExitTwo}, for which we make use of the
following auxiliary lemma:

\blems
\label{LemAuxiliaryTwo}
For $u \in \dom \gfun \cap \real^+$, we have the decrement bound
\begin{align}
\label{EqnDampedDecrementBound}
\gfun(\ustep) \le \gfun(0) + \ustep \inprod{\nabla f(x)}{\vsketch} -
\ustep \| [\nabla^2 f(x)]^{1/2} \vsketch \|_2 - \log \big(1- \ustep \|
       [\nabla^2 f(x)]^{1/2} \vsketch \|_2 \big).
\end{align}
provided that $u\|[\nabla^2 f(x)]^{1/2} \vsketch \|_2<1$.
\elems
\blems
\label{LemAuxiliaryThree}
With probability at least $1 - c_1 e^{-c_2 \numproj}$, we have
\begin{align}
\label{EqnAuxiliaryTwo} 
\| [\nabla^2 f(x)]^{1/2} \vsketch \|_2^2 & \leq
\left(\frac{1+\epsilon}{1-\epsilon} \right)^2 \,\,
\big[\newdecsketch_f(x) \big]^2 \,.
\end{align}
\elems

\noindent The proof of these lemmas are provided in
Appendices~\ref{AppLemAuxiliaryTwo} and~\ref{AppLemAuxiliaryThree}.
Using them, let us prove the claims~\eqref{EqnDampedDecrementExit}
and~\eqref{EqnDampedDecrementExitTwo}.  Recalling our shorthand
$\epsilon^\prime \defn \frac{1+\epsilon}{1-\epsilon}-1 =
\frac{2\epsilon}{1-\epsilon}$, substituting
inequality~\eqref{EqnAuxiliaryTwo} into the decrement
formula~\eqref{EqnDampedDecrementBound} yields
\begin{align}
\gfun(\ustep) & \leq \gfun(0) - \ustep \newdecsketch_f(x)^2 -
\ustep (1+\epsilon^\prime) \,\,\newdecsketch_f(x) - \log(1 -
\ustep (1+\epsilon^\prime)
\,\,\newdecsketch_f(x)) \label{EqnApproxNewtonStepExpansion} \\ 
& = \gfun(0) - {\left\{ \ustep (1+\epsilon^\prime)^2
  \newdecsketch_f(x)^2 + \ustep (1+\epsilon^\prime)
  \,\,\newdecsketch_f(x) + \log(1 - \ustep (1+\epsilon^\prime)
  \,\,\newdecsketch_f(x)) \right\}} \nonumber \\ & \qquad \qquad
\qquad \qquad + \ustep(({1+\epsilon^\prime})^2 -1
)\newdecsketch_f(x)^2 \nonumber\,
\end{align}
where we added and subtracted $\ustep (1+\epsilon^\prime)^2
\newdecsketch_f(x)^2$ so as to obtain the final equality.

We now prove inequality~\eqref{EqnDampedDecrementExit}.  Now setting
$\ustep = \shat \defn \frac{1}{1+(1+\epsilon^\prime)
  \newdecsketch_f(x)}$, which satisfies the conditions of Lemma
\ref{LemAuxiliaryTwo} yields
\begin{align*}
\gfun(\shat) \le \gfun(0) - (1+\epsilon^\prime) \,\,\newdecsketch_f(x)
+ \log(1+ (1+\epsilon^\prime) \,\,\newdecsketch_f(x)) +
\frac{({\epsilon^\prime}^2 + 2 \epsilon^\prime
  )\newdecsketch_f(x)^2}{1+(1+\epsilon^\prime) \newdecsketch_f(x) }\,.
\end{align*}
Making use of the standard inequality $-u + \log(1+u)\le -
\frac{\frac{1}{2}u^2}{(1+u)}$ (for instance, see the
book~\cite{Boyd02}), we find that
\begin{align*}
\gfun(\shat) & \leq \gfun(0) - \frac{\frac{1}{2}
  (1+\epsilon^\prime)^2\newdecsketch_f(x)^2}{1+(1+\epsilon^\prime)
  \newdecsketch_f(x) } + \frac{({\epsilon^\prime}^2 + 2
  \epsilon^\prime )\newdecsketch_f(x)^2}{1+(1+\epsilon^\prime)
  \newdecsketch_f(x) } \\ & = \gfun(0) -
(\frac{1}{2}-\frac{1}{2}{\epsilon^\prime}^2 - \epsilon^\prime
)\newdecsketch_f(x)^2 \shat \\
& \leq \gfun(0) - \alpha \newdecsketch_f(x)^2 \shat,
\end{align*}
where the final inequality follows from our assumption $\alpha\le
\frac{1}{2}- \frac{1}{2}{\epsilon^\prime}^2 - \epsilon^\prime$. This
completes the proof of the bound~ \eqref{EqnDampedDecrementExit}.
Finally, the lower bound~\eqref{EqnDampedDecrementExitTwo} follows by
setting $\ustep = \lineb \shat$ into the decrement
inequality~\eqref{EqnDampedDecrementBound}. \\


\subsubsection{Proof of Lemma~\ref{LemKeyDecrement}(b)}

The proof of this part hinges on the following auxiliary lemma:
\blems
\label{LemBeanery}
For all $\epsilon \in (0, 1/2)$, we have
\begin{subequations}
\begin{align} 
\label{EqnNewtonLinQuadCond}
   \newdec_f(\xsketch) & \leq \frac{(1 + \epsilon) \newdec^2_f(x) +
     \epsilon \newdec_f(x)}{\Big (1 -( 1 + \epsilon) \newdec_f(x)
     \Big)^2}, \qquad \mbox{and} \\
 \label{EqnNewtonLinQuadCond2} 
(1 - \epsilon) \, \newdec_f(\xsketch) \; \leq \;
 \newdecsketch_f(\xsketch) & \leq (1+\epsilon) \newdec_f(\xsketch) \,,
\end{align}
\end{subequations}
where all bounds hold with probability at least $1-c_1 e^{-c_2
  \numproj \epsilon^2}$.
\elems
\noindent See Appendix~\ref{AppLemBeanery} for the proof. \\

\noindent We now use Lemma~\ref{LemBeanery} to prove the two claims
in the lemma statement.

\paragraph{Proof of the bound~\eqref{EqnSecondPhaseControl}:}
 Recall from the theorem statement that $\NEWETA \defn \frac{1}{8}
 \,\frac{1-\frac{1}{2}(\frac{1 + \epsilon}{1-\epsilon})^2 -
   \linea}{(\frac{1+\epsilon}{1 - \epsilon})^3}$.  By examining the
 roots of a polynomial in $\epsilon$, it can be seen that $\NEWETA
 \leq \frac{1-\epsilon}{1 + \epsilon} \, \frac{1}{16}$.
\begin{align*}
(1 + \epsilon) \, \newdec_f(\xit{t}) \leq (1 + \epsilon)
  \newdecsketch_f(\xit{t}) \leq (1 + \epsilon) \, \NEWETA \leq
  \frac{1}{16}
\end{align*}
By applying the inequalities~\eqref{EqnNewtonLinQuadCond2}, we have
\begin{align}
(1 + \epsilon) \newdec_f(x) \leq \frac{1 + \epsilon}{1-\epsilon}
  \newdecsketch_f(x) \; \leq 
\frac{1 + \epsilon}{1-\epsilon} \, \NEWETA \, \leq \, \frac{1}{16}
\end{align}
 whence inequality~\eqref{EqnNewtonLinQuadCond} implies that
\begin{align}
\label{EqnNewtonDampedLinearConvergence}
 \newdec_f(\xsketch) & \leq \frac{ \frac{1}{16} \newdec_f(x) +
   \epsilon \newdec_f(x)}{(1 - \frac{1}{16})^2} \; \leq \;
 \left(\frac{16}{225} + \frac{256}{225}\epsilon \right )\newdec_f(x) \; \leq \;
 \frac{16}{25}  \newdec_f(x).
\end{align}
Here the final inequality holds for all $\epsilon \in (0, 1/2)$.
Combining the bound~\eqref{EqnNewtonLinQuadCond2} with
inequality~\eqref{EqnNewtonDampedLinearConvergence} yields
\begin{align*}
\newdecsketch_f(\xsketch) \leq (1+\epsilon) \newdec_f(\xsketch) & \leq
(1 + \epsilon) \big( \frac{16}{25} \big) \; \newdecsketch_f(x) \; \leq
\; \newdecsketch_f(x)\,,
\end{align*}
where the final inequality again uses the condition $\epsilon \in
(0,\frac{1}{2})$.  This completes the proof of the
bound~\eqref{EqnSecondPhaseControl}.


\paragraph{Proof of the bound~\eqref{EqnSecondPhaseContract}:}
This inequality has been established as a consequence of proving the
bound~\eqref{EqnNewtonDampedLinearConvergence}.


\subsection{Proof of Theorem~\ref{ThmSelfConcordantBarrier}}

Given the proof of Theorem~\ref{ThmSelfConcordantUnconstrained}, it
remains only to prove the following modified version of
Lemma~\ref{LemApproxNewtonStep}.  It applies to the exact and sketched
Newton directions $\vnewt, \vsketch \in \real^\usedim$ that are
defined as follows
\begin{subequations}
\begin{align}
\label{EqnNewtonOrgModified} 
\vnewt & \defn \arg \min_{z \in \Constraint-x} \ENCMIN{ \frac{1}{2}\|
  \fpp(x)^{1/2} z\|_2^2 + \inprod{z}{\nabla f(x)} + \frac{1}{2}
  \inprod{z}{ \nabla^2 g(x) z} }, \\
\label{EqnNewtonApproxModified}
\vsketch & = \arg\min_{z \in \Constraint-x} \ENCMIN{
  \underbrace{\frac{1}{2}\| \Sketch \fpp(x)^{1/2} z\|_2^2 +
    \inprod{z}{\nabla f(x)}+ \frac{1}{2} \inprod{z} {\nabla^2 g(x)
      z}}_{\Psi(z; \Sketch)} }.
\end{align}
\end{subequations}
Thus, the only difference is that the Hessian $\nabla^2 f(x)$ is
sketched, whereas the term $\nabla^2 g(x)$ remains unsketched.

\blems 
\label{LemApproxNewtonStepModified}
Let $\Sketch \in \real^{\numproj \times \numobs}$ be a sub-Gaussian or
ROS sketching matrix, and let $x \in \real^\usedim$ be a (possibly
random) vector independent of $\Sketch$.  If $\numproj \ge c_0
\max_{x\in \Constraint}\frac{\Width(\fpp(x)^{1/2}
  \ConeSet)^2}{\epsilon^2}$, then
\begin{align} 
\left\| \fpp(x)^{1/2} (\vsketch - \vnewt) \right\|_2 & \leq \epsilon
\left\| \fpp(x)^{1/2} \vnewt \right\|_2~
\end{align}
with probability at least $1-c_1e^{-c_2 \numproj \epsilon^2}$.
\elems
%

\section{Discussion}
\label{SecDiscussion}

In this paper, we introduced and analyzed the Newton sketch, a
randomized approximation to the classical Newton updates.  This
algorithm is a natural generalization of the Iterative Hessian Sketch
(IHS) updates analyzed in our earlier work~\cite{PilWai14b}.  The IHS
applies only to constrained least-squares problems (for which the
Hessian is independent of the iteration number), whereas the Newton
Sketch applies to any any twice differentiable function subject to a
closed convex constraint set.  We described various applications of
the Newton sketch, including its use with barrier methods to solve
various forms of constrained problems.  For the minimization of
self-concordant functions, the combination of the Newton sketch within
interior point updates leads to much faster algorithms for an
extensive body of convex optimization problems.

Each iteration of the Newton sketch always has lower computational
complexity than classical Newton's method.  Moreover, it has lower
computational complexity than first-order methods when either $\numobs
\ge \usedim^2$ or $\usedim \geq \numobs^2$ (using the dual strategy);
here $\numobs$ and $\usedim$ denote the dimensions of the data matrix
$\Amat$. In the context of barrier methods, the parameters $\numobs$
and $\usedim$ typically correspond to the number of constraints and
number of variables, respectively. In many ``big data'' problems, one
of the dimensions is much larger than the other, in which case the
Newton sketch is advantageous.  Moreover, sketches based on the
randomized Hadamard transform are well-suited to in parallel
environments: in this case, the sketching step can be done in
$\order(\log \numproj)$ time with $\order (\numobs \usedim)$
processors. This scheme significantly decreases the amount of central
computation---namely, from $\order(m^2 d + nd \log m)$ to
$\order(\numproj^2 \usedim + \log \usedim)$.

There are a number of open problems associated with the Newton sketch.
Here we focused our analysis on the cases of sub-Gaussian and
randomized orthogonal system (ROS) sketches.  It would also be
interesting to analyze sketches based on coordinate sampling, or other
forms of ``sparse'' sketches (for instance, see the
paper~\cite{kane2014sparser}).  Such techniques might lead to
significant gains in cases where the data matrix $\Amat$ is itself
sparse: more specifically, it may be possible to obtain sketched
optimization algorithms whose computational complexity only scales
with number of nonzero entries in the data matrices the full
dimensionality $\numobs \usedim$.  Finally, it would be interesting
to explore the problem of lower bounds on the sketch dimension
$\numproj$.  In particular, is there a threshold below which any
algorithm that has access only to gradients and $\numproj$-sketched
Hessians must necessarily converge at a sub-linear rate, or in a way
that depends on the strong convexity and smoothness parameters?  Such
a result would clarify whether or not the guarantees in this paper are
improvable.



\subsection*{Acknowledgements}
Both authors were partially supported by Office of Naval Research MURI
grant N00014-11-1-0688, and National Science Foundation Grants
CIF-31712-23800 and DMS-1107000. In addition, MP was supported by a
Microsoft Research Fellowship.


\appendix

\section{Technical results for Theorem~\ref{ThmSelfConcordantUnconstrained}}

In this appendix, we collect together various technical results and
proofs that are required in the proof of
Theorem~\ref{ThmSelfConcordantUnconstrained}.

\subsection{Proof of Lemma~\ref{LemApproxNewtonStep}}

Let $u$ be a unit-norm vector independent of $\Sketch$, and consider
the random quantities
\begin{subequations}
\begin{align}
\label{EqnDefnZinf}
\ZINF(\Sketch,x) & \defn \inf_{v \in \fpp(x)^{1/2}\KCONE^t \cap
  \Sphere{\numobs}} \| \Sketch v\|_2^2 \quad \mbox{and} \\
\label{EqnDefnZsup}
\ZSUP(\Sketch,x) & \defn \sup_{v \in \fpp(x)^{1/2}\KCONE^t \cap \Sphere{\numobs}} \Big|
\inprod{\fixvec}{(\Sketch^T \Sketch - I_\numobs) \,
  v} \Big|.
\end{align}
\end{subequations}
By the optimality and feasibility of $\vsketch$ and $\vnewt$
(respectively) for the sketched Newton update~\eqref{EqnNewtonApprox},
we have
\begin{align*}
\frac{1}{2} \|\Sketch \fpp(x)^{1/2} \vsketch\|_2^2 -
\inprod{\vsketch}{\fp(x)} & \leq \frac{1}{2} \|\fpp(x)^{1/2}
\vnewt\|_2^2 - \inprod{\vnewt}{\fp(x)}.
\end{align*}

Defining the difference vector $\errvec \defn \vsketch - \vnewt$, some
algebra leads to the basic inequality
\begin{align}
\label{EqnBasicInequality}
\frac{1}{2} \|\Sketch \fpp(x)^{1/2} \errvec\|_2^2 & \leq -
\inprod{\fpp(x)^{1/2} \vnewt}{\Sketch^T \Sketch \fpp(x)^{1/2} \errvec}
+ \inprod{\errvec}{\fp(x)}.
\end{align}
Moreover, by the optimality and feasibility of $\vnewt$ and $\vsketch$
for the exact Newton update~\eqref{EqnNewtonOrig}, we have
\begin{align}
\label{EqnOptimalForOriginal}
\inprod{\fpp(x) \vnewt - \fp(x) }{\errvec} \; = \; \inprod{\fpp(x)
  \vnewt - \fp(x) }{\vsketch - \vnewt} & \geq 0.
\end{align}
Consequently, by adding and subtracting $\inprod{\fpp(x)
  \vnewt}{\errvec}$, we find that
\begin{align}
\label{EqnModifyForFixedSketchConcordant}
\frac{1}{2 } \|\Sketch \fpp(x)^{1/2} \errvec\|_2^2 & \leq \Big|
\inprod{\fpp(x)^{1/2} \vnewt} {\big (I_\numobs - \Sketch^T \Sketch
  \big) \fpp(x)^{1/2} \errvec }\Big|\,.
\end{align}
By definition, the error vector $\errvec$ belongs to the cone $\KCONE^t$
and the vector $\fpp(x)^{1/2} \vnewt$ is fixed and independent of the
sketch.  Consequently, invoking definitions~\eqref{EqnDefnZinf}
and~\eqref{EqnDefnZsup} of the random variables $\ZINF$ and $\ZSUP$
yields
\begin{align*}
\frac{1}{2} \|\Sketch \fpp(x)^{1/2} \errvec\|_2^2 &\geq
\frac{\ZINF}{2} \| \fpp(x)^{1/2} \errvec\|_2^2, \\
~ \Big| \inprod{\fpp(x)^{1/2} \vnewt} {\big (I_\numobs - \Sketch^T
  \Sketch \big) \fpp(x)^{1/2} \errvec} \Big| &\leq \ZSUP
\|\fpp(x)^{1/2} \vnewt\|_2 \, \|\fpp(x)^{1/2} \errvec\|_2,
\end{align*}
Putting together the pieces, we find that
\begin{align} 
\label{EqnLemSketchZ} 
\left\| \fpp(x)^{1/2} (\vsketch-\vnewt) \right\|_2 \le \frac{2
  Z_2(\Sketch,x)}{Z_1(\Sketch,x)}\left\| \fpp(x)^{1/2} (\vnewt)
\right\|_2~.
\end{align}
Finally, for any $\delta \in (0,1)$, let us define the event
$\Event(\delta) = \{\ZINF \geq 1-\delta, \quad \mbox{and} \quad \ZSUP
\leq \delta \}$.  By Lemma 4 and Lemma 5 from our previous
paper~\cite{PilWai14a}, we are guaranteed that $\mprob[\Event(\delta)]
\geq 1 - c_1 e^{-c_2 \numproj \delta^2}$.  Conditioned on the event
$\Event(\delta)$, the bound~\eqref{EqnLemSketchZ} implies that
\begin{align*}
\left\| \fpp(x)^{1/2} (\vsketch-\vnewt) \right\|_2 \le \frac{2
  \delta}{1-\delta}\left\| \fpp(x)^{1/2} (\vnewt) \right\|_2~.
\end{align*}
By setting $\delta = \frac{\epsilon}{4}$, the claim follows.


\subsection{Proof of Lemma~\ref{LemAuxiliaryTwo}}
\label{AppLemAuxiliaryTwo}

By construction, the function $\gfun(u) = f(x + u \vsketch)$ is
strictly convex and self-concordant.  Consequently, it satisfies the
bound $\frac{d}{d u} \left(\gfun^\dprime(u)^{-1/2}\right) \leq 1$,
whence
\begin{align*}
\gfun^\dprime(s)^{-1/2} - \gfun^\dprime(0)^{-1/2} \; = \; \int_{0}^s
\frac{d}{du} \left( \gfun^\dprime(u)^{-1/2} \right) du \; \leq \; s.
\end{align*}
or equivalently $\gfun^\dprime(s) \le \frac{\gfun^\dprime(0)}{(1-s
  \gfun^\dprime(0)^{1/2})^2}$ for $s \in \dom \gfun \cap
[0,g^\dprime(0)^{-1/2})$. Integrating this inequality twice yields the
  bound
\begin{align} 
\label{EqnNewtonStepExpansion}
\gfun(u) \le \gfun(0) + u \gfun^\prime(0) - u \gfun^\dprime(0)^{1/2} -
\log(1-u \gfun^\dprime(0)^{1/2})\,.
\end{align}
Since $\gfun^\prime(u) = \inprod{\nabla f(x + u \vsketch)}{\vsketch}$
and $\gfun^{\dprime}(u) = \inprod{\vsketch}{\nabla^2 f(x + u\vsketch)
  \vsketch}$, the decrement bound~\eqref{EqnDampedDecrementBound}
follows.\\


\subsection{Proof of Lemma~\ref{LemAuxiliaryThree}}
\label{AppLemAuxiliaryThree}

We perform this analysis conditional on the bound~\eqref{EqnLemSketch}
from Lemma~\ref{LemApproxNewtonStep}.  We begin by observing that
\begin{align}
\| [\nabla^2 f(x)]^{1/2} \vsketch \|_2 & \leq \| [\nabla^2 f(x)]^{1/2}
\vnewt \|_2 + \| [\nabla^2 f(x)]^{1/2} (\vsketch-\vnewt) \|_2
\nonumber \\
\label{EqnDampedDecrementBoundSecondTerm}
& = \newdec_f(x) + \| [\nabla^2 f(x)]^{1/2} (\vsketch-\vnewt )\|_2\,.
\end{align}
Lemma~\ref{LemApproxNewtonStep} implies that $\|\nabla^2
[f(x)]^{1/2}(\vsketch-\vnewt)\|_2 \leq \epsilon \|\nabla^2
[f(x)]^{1/2} \vnewt \|_2 = \epsilon \newdec_f(x)$.
In conjunction with the
bound~\eqref{EqnDampedDecrementBoundSecondTerm}, we see that
\begin{align} 
\label{EqnDampedDecrementBoundSecondTerm}
\| [\nabla^2 f(x)]^{1/2} \vsketch \|_2 & \leq (1+\epsilon)
\newdec_f(x)\,.
\end{align}
Our next step is to lower bound the term $\inprod{\nabla
  f(x)}{\vsketch}$: in particular, by adding and subtracting a factor
of the original Newton step $\vnewt$, we find that
\begin{align} 
\inprod{\nabla f(x)}{\vsketch} &= \inprod{[\nabla^2 f(x)]^{-1/2}
  \nabla f(x)}{\nabla^2 [f(x)]^{1/2}\vsketch} \nonumber \\ 
& = \inprod{[\nabla^2 f(x)]^{-1/2}\nabla f(x)}{\nabla^2
  [f(x)]^{1/2}\vnewt} + \inprod{[\nabla^2 f(x)]^{-1/2}\nabla
  f(x)}{\nabla^2 [f(x)]^{1/2}(\vsketch-\vnew)} \nonumber \\
& = -\| \nabla^2 [f(x)]^{-1/2} \nabla f(x) \|_2^2 + \inprod{[\nabla^2
    f(x)]^{-1/2} \nabla f(x)}{\nabla^2 [f(x)]^{1/2} (\vsketch -
  \vnewt)} \nonumber \\
& \leq -\| \nabla^2 [f(x)]^{-1/2}\nabla f(x) \|_2^2 + {\|[\nabla^2
    f(x)]^{-1/2} \nabla f(x)\|_2} \|\nabla^2
          [f(x)]^{1/2}(\vsketch-\vnewt)\|_2 \nonumber \\
& = -\newdec_f(x)^2 + \newdec_f(x) \|\nabla^2
      [f(x)]^{1/2}(\vsketch-\vnewt)\|_2 \nonumber \\ 
& \leq -\newdec_f(x)^2(1-\epsilon),
\end{align}
where the final step again makes use of
Lemma~\ref{LemApproxNewtonStep}.  Repeating the above argument in the
reverse direction yields the lower bound $\inprod{\nabla
  f(x)}{\vsketch} \geq -\newdec_f(x)^2 (1+\epsilon)$, so that we may
conclude that
\begin{align}
\label{EqnHeadphones}
\vert \newdecsketch_f(\xit{t}) - \newdec_f(\xit{t}) \vert \le \epsilon
\newdec_f(\xit{t}).
\end{align}
Finally, by squaring both sides of the
inequality~\eqref{EqnDampedDecrementBoundSecondTerm} and combining
with the above bounds gives
\begin{align*}
\| [\nabla^2 f(x)]^{1/2} \vsketch \|_2^2 \; \leq \;
\frac{-\,(1+\epsilon)^2}{1-\epsilon} \, \inprod{\nabla f(x)}{\vsketch}
& = \frac{(1+\epsilon)^2}{1-\epsilon} \,\,\newdecsketch^2_f(x) \; \leq
\; \left(\frac{1+\epsilon}{1-\epsilon} \right)^2 \,\,
\newdecsketch^2_f(x),
\end{align*}
as claimed.


\subsection{Proof of Lemma~\ref{LemBeanery}}
\label{AppLemBeanery}

We have already proved the bound~\eqref{EqnNewtonLinQuadCond2} during
our proof of Lemma~\ref{LemAuxiliaryThree}---in particular, see
equation~\eqref{EqnHeadphones}.  Accordingly, it remains only to prove
the inequality~\eqref{EqnNewtonLinQuadCond}.

Introducing the shorthand $\SHORTA \defn (1+\epsilon)\lambda_f(x)$, we
first claim that the Hessian satisfies the sandwich relation
\begin{align}
\label{EqnHessianBound}
(1- \stepsize \alpha)^2 \nabla^2 f(x) \preceq \nabla^2 f(x+ \stepsize
\vsketch) \preceq \frac{1}{(1- \stepsize\alpha)^2} \nabla^2 f(x)\,,
\end{align}
for $\vert 1-s\alpha \vert<1$ where $\alpha = (1+\epsilon)\lambda_f(x)$, with probability at least $1-c_1 e^{-c_2 \numproj \epsilon^2}$. Let us recall
Theorem 4.1.6 of Nesterov~\cite{Nesterov04}: it guarantees that
\begin{align}
\label{EqnHessianBoundNesterov}
(1- \stepsize \|\vsketch\|_x)^2 \nabla^2 f(x) \preceq \nabla^2 f(x +
\stepsize \vsketch) \preceq \frac{1}{(1 - \stepsize \|\vsketch\|_x)^2}
\nabla^2 f(x)~.
\end{align}
Now recall the bound~\eqref{EqnLemSketch} from
Lemma~\ref{LemApproxNewtonStep}: combining it with an application of
the triangle inequality (in terms of the semi-norm \mbox{$\|v\|_x = \|
  \fpp(x)^{1/2} v \|_2$)} yields
\begin{align*} 
\left\| \fpp(x)^{1/2} \vsketch \right\|_2 \le & (1+\epsilon) \left\|
\fpp(x)^{1/2} \vnewt \right\|_2~ \; = \; (1+\epsilon) \|\vnewt \|_x\,,
\end{align*}
with probability at least $1-e^{-c_1 \numproj \epsilon^2}$, and
substituting this inequality into the
bound~\eqref{EqnHessianBoundNesterov} yields the sandwich
relation~\eqref{EqnHessianBound} for the Hessian. \\

Using this sandwich relation~\eqref{EqnHessianBound}, the Newton
decrement can be bounded as
\begin{align*}
\lambda_f(\xsketch) & = \| \nabla^2f(\xsketch)^{-1/2} \nabla f(\xsketch)\|_2 \\
& \leq \frac{1}{\left(1-(1+\epsilon)\lambda_f(x)\right)} \| \nabla^2
f(x)^{-1/2} \nabla f(\xsketch)\|_2 \\
& = \frac{1}{\left(1-(1+\epsilon)\lambda_f(x)\right)} \left\| \nabla^2
f(x)^{-1/2} \left( \nabla f(x)+\int_0^1 \nabla^2 f(x + \stepsize
\vsketch)\vsketch \; d \stepsize \right)\right\|_2 \\
& = \frac{1}{\left(1-(1+\epsilon)\lambda_f(x)\right)} \left\| \nabla^2
f(x)^{-1/2}\left(\nabla f(x)+\int_0^1 \nabla^2 f(x + \stepsize
\vsketch) \vnewt \: d \stepsize + \Delta \right) \right\|_2,
\end{align*}
where we have defined $\Delta = \int_0^1 \fpp(x + \stepsize \vsketch)
\, (\vsketch-\vnewt) \, d \stepsize$.  By the triangle inequality, we
can write $\lambda_f(\xsketch) \leq
\frac{1}{\left(1-(1+\epsilon)\lambda_f(x)\right)} \big( M_1 + M_2
\big)$, where
\begin{align*}
M_1 & \defn \left\| \nabla^2 f(x)^{-1/2} \left ( \nabla f(x)+\int_0^1
\nabla^2 f(x+t \vsketch) \vnewt dt \right) \right\|_2, \quad
\mbox{and} \quad M_2 \defn \left \|\fpp(x)^{-1/2} \Delta \right \|_2.
\end{align*}
In order to complete the proof, it suffices to show that
\begin{align*}
M_1 & \leq
\frac{(1+\epsilon)\lambda_f(x)^2}{1-(1+\epsilon)\lambda_f(x)}, \quad
\mbox{and} \quad M_2 \leq \frac{\epsilon
  \newdec_f(x)}{1-(1+\epsilon)\lambda_f(x)}.
\end{align*}

\paragraph{Bound on $M_1$:}
Re-arranging and then invoking the Hessian sandwich
relation~\eqref{EqnHessianBound} yields
\begin{align*}
M_1 & = \left\| \int_0^1 \left ( \fpp(x)^{-1/2} \nabla^2 f(x
+\stepsize \vsketch) \fpp(x)^{-1/2} -I \right) d \stepsize \;
\left(\nabla^2 f(x)^{1/2} \vnewt \right) \right \|_2 \\
& \leq \left | \int_0^1 \left(\frac{1}{(1 - \stepsize (1 + \epsilon)
  \lambda_f(x))^2}-1 \right) d \stepsize \; \right | \; \left |
\left(\nabla^2 f(x)^{1/2} \vnewt \right) \right \|_2 \\
& = \frac{(1 + \epsilon) \lambda_f(x)}{ 1 - (1 + \epsilon)
  \lambda_f(x)} \left\| \nabla^2 f(x)^{1/2} \vnewt \right \|_2 \\
& = \frac{(1+\epsilon)\lambda^2_f(x)}{1-(1+\epsilon)\lambda_f(x)}.
\end{align*}

\paragraph{Bound on $M_2$:}
We have
\begin{align*}
M_2 & = \left \| \int_0^1 \fpp(x)^{-1/2}\fpp(x + \stepsize \vsketch)
\fpp(x)^{-1/2} d \stepsize \fpp(x)^{1/2} (\vsketch - \vnewt)
\right\|_2 \\
& \le \left\| \int_0^1 \frac{1}{(1 - \stepsize
  (1+\epsilon)\lambda_f(x))^2} d \stepsize \fpp(x)^{1/2} (\vsketch-
\vnewt) \right\|_2 \nonumber \\
& = \frac{1}{1-(1+\epsilon)\newdec_f(x)} \left\| \fpp(x)^{1/2}
(\vsketch - \vnewt) \right\|_2 \nonumber \\
 & \stackrel{(i)}{\leq} \frac{1}{1-(1+\epsilon)\newdec_f(x)} \epsilon
\left\| \fpp(x)^{1/2} \vnewt \right\|_2 \\
& = \frac{\epsilon \newdec_f(x)}{1-(1+\epsilon)\lambda_f(x)} \nonumber
\,,
\end{align*}
where the inequality in step (i) follows from
Lemma~\ref{LemApproxNewtonStep}.
%


\section{Proof of Lemma~\ref{LemApproxNewtonStepModified}}

The proof follows the basic inequality argument of the proof of
Lemma~\ref{LemApproxNewtonStep}.  Since $\vsketch$ and $\vnewt$ are
optimal and feasible (respectively) for the sketched Newton
problem~\eqref{EqnNewtonApproxModified}, we have $\Psi(\vsketch;
\Sketch) \leq \Psi(\vnewt; \Sketch)$.  Defining the difference vector
$\errvec \defn \vsketch - v$, some algebra leads to the basic
inequality
\begin{align*}
\frac{1}{2} \|\Sketch \fpp(x)^{1/2} \errvec\|_2^2 + \frac{1}{2}
\inprod{\errvec}{\nabla^2 g(x) \errvec} & \leq - \inprod{\fpp(x)^{1/2}
  \vnewt}{\Sketch^T \Sketch \fpp(x)^{1/2} \errvec} \\
& \qquad \qquad \qquad \qquad + \inprod{\errvec}{\big (
  \fp(x)-\nabla^2 g(x) \big) \vnewt}.
\end{align*}
On the other hand since $\vnewt$ and $\vsketch$ are optimal and
feasible (respectively) for the Newton
step~\eqref{EqnNewtonOrgModified}, we have
\begin{align*}
\inprod{\fpp(x) \vnewt + \nabla^2 g(x) \vnewt - \fp(x) }{\errvec} & \geq 0.
\end{align*}
Consequently, by adding and subtracting $\inprod{\fpp(x)
  \vnewt}{\errvec}$, we find that
\begin{align*}
\frac{1}{2 } \|\Sketch \fpp(x)^{1/2} \errvec\|_2^2 + \frac{1}{2}
\inprod{\vnewt}{\nabla^2 g(x) \vnewt} & \leq \Big|
\inprod{\fpp(x)^{1/2} \vnewt} {\big (I_\numobs - \Sketch^T \Sketch
  \big) \fpp(x)^{1/2} \errvec }\Big|\,.
\end{align*}
By convexity of $g$, we have $\nabla^2 g(x) \succeq 0$, whence
\begin{align*}
\frac{1}{2 } \|\Sketch \fpp(x)^{1/2} \errvec\|_2^2 & \leq \Big|
\inprod{\fpp(x)^{1/2} \vnewt} {\big (I_\numobs - \Sketch^T \Sketch
  \big) \fpp(x)^{1/2} \errvec }\Big|\,
\end{align*}
Given this inequality, the remainder of the proof follows as in the
proof of Lemma~\ref{LemApproxNewtonStep}.
%

\section{Gaussian widths with $\ell_1$-constraints}
\label{AppEllOneWidth}

In this appendix, we state and prove an elementary lemma that bounds
for the Gaussian width for a broad class of $\ell_1$-constrained
problems.  In particular, given a twice-differentiable convex function
$\LinkFun$, a vector $c \in \real^\usedim$, a radius $R$ and a
collection of $\usedim$-vectors $\{a_i\}_{i=1}^\numobs$, consider a
convex program of the form
\begin{align}
\label{EqnEllOneGLM}
\min_{ \substack{x \in \Constraint}} \Big \{ \sum_{i=1}^\numobs
\LinkFun \big( \inprod{a_i}{x} \big) + \inprod{c}{x} \Big \}, \qquad
\mbox{where} \quad \Constraint = \{x \in \real^\usedim \, \mid \,
\|x\|_1 \leq R \}.
\end{align}

\blems
\label{LemL1ConstrainedGLMWidth}
Suppose that the $\ell_1$-constrained program~\eqref{EqnEllOneGLM} has
a unique optimal solution $\xstar$ such that $\|\xstar\|_0\le s$ for
some integer $\kdim$.  Then denoting the tangent cone at $\xstar$ by
$\ConeSet$, then
\begin{align*}
\max_{x \in \Constraint} \Width(\nabla^2 f(x)^{1/2} \ConeSet) & \leq 6
\sqrt{\kdim \log \usedim} \;
\sqrt{\frac{\LinkFunDouble_{\max}}{\LinkFunDouble_{\min}}} \: \;
\frac{\max \limits_{j=1,\ldots, \usedim} \|A_j\|_2
}{\sqrt{\sigkminsq(A)}}\,,
\end{align*}
where 
\begin{align*}
\LinkFunDouble_{\min} = \min \limits_{x \in \Constraint} \min
\limits_{ i = 1, \ldots, \numobs} \, \LinkFunDouble(\inprod{a_i}{x},
y_i), \quad \mbox{and} \quad \LinkFunDouble_{\max} = \max \limits_{x
  \in \Constraint} \; \max \limits_{i = 1, \ldots, \numobs} \,
\LinkFunDouble(\inprod{a_i}{x}, y_i).
\end{align*}

\elems 
\proof 
It is well known~\cite{NegRavWaiYu12,PilWai14a} that the tangent cone
of the $\ell_1$-norm at any $\kdim$-sparse solution is a subset of the
cone $\{z \in \real^\usedim \, \mid \, \|z\|_1 \leq 2 \sqrt{\kdim}
\|z\|_2\}$.  Using this fact, we have the following sequence of upper
bounds
\begin{align*}
 \Width(\nabla^2 f(x)^{1/2} \ConeSet) &= \Exs_w \max_{\substack{z^T
     \nabla^2 f(x) z = 1\,,\\ z \in \ConeSet}} \, \inprod{w}{ \nabla^2
   f(x)^{1/2} z} \\
& = \Exs_w \max_{\substack{z^T A^T \diag\left(\LinkFunDouble(
     \inprod{a_i}{x} x,y_i)\right) \Amat z = 1\,, \\ z \in \ConeSet}}
 \, \inprod{w}{ \diag \left( \LinkFunDouble(\inprod{a_i}{x},
   y_i)\right)^{1/2} \Amat z} \\
& \leq \Exs_w \max_{\substack{z^T A^T A z \le 1/\LinkFunDouble_{\min}
     \\ z \in \ConeSet}} \, \inprod{w}{ \diag \left(
   \LinkFunDouble(\inprod{a_i}{x}, y_i) \right)^{1/2} \Amat z} \\
& \leq \Exs_w \max_{\|z\|_1 \leq
   \frac{2\sqrt{\kdim}}{\sqrt{\sigkminsq(\Amat)}}\frac{1}{\sqrt{\LinkFunDouble_{\min}}}}
 \, \inprod{w}{ \diag\left(
   \LinkFunDouble(\inprod{a_i}{x},y_i)\right)^{1/2}Az} \\
& = \frac{2\sqrt{\kdim}}{\sqrt{\sigkminsq(A)}}
 \frac{1}{\sqrt{\LinkFunDouble_{\min}}} \Exs_w \, \|\Amat^T
 \diag\left(\LinkFunDouble( \inprod{a_i}{x}, y_i) \right )^{1/2}
 w\|_{\infty} \\
& =
 \frac{2\sqrt{s}}{\sqrt{\sigkminsq(A)}}\frac{1}{\sqrt{\LinkFunDouble_{\min}}}
 \Exs_w \, \max_{j = 1, \ldots, \usedim}\, \Big\vert
 \underbrace{\sum_{i = 1, \ldots, \numobs} w_i A_{ij} \LinkFunDouble(
   \inprod{a_i}{x}, y_i)^{1/2} }_{Q_{j}} \Big \vert.
\end{align*}
Here the random variables $Q_{j}$ are zero-mean Gaussians with
variance at most
\begin{align*}\sum_{i=1, \ldots, \numobs} A_{ij}^2 \LinkFunDouble(
\inprod{a_i}{x}, y_i) \leq \LinkFunDouble_{\max}
\|A_j\|_2^2. 
\end{align*}
Consequently, applying standard bounds on the suprema of Gaussian
variates~\cite{LedTal91}, we obtain
\begin{align*}
\Exs_w \, \max_{j=1,\ldots, \usedim}\, \Big \vert {\sum_{i=1,\ldots,
    \numobs} w_i A_{ij} \LinkFunDouble( \inprod{a_i}{x}, y_i)^{1/2} }
\Big \vert \leq 3 \sqrt{\log \usedim} \sqrt{{ \LinkFunDouble_{\max}}}
     {\max_{j=1, \ldots,\usedim} \|A_j\|_2 }.
\end{align*}
When combined with the previous inequality, the claim follows.
\endproof



\bibliographystyle{plain}

\bibliography{mert_super}
\end{document}